\documentclass[12pt]{article}
\usepackage{amsfonts, amsmath, amssymb, amsgen, amsthm, amscd,
latexsym}
\usepackage[all]{xy}

\def\Diff{\mathop{\rm Diff}\nolimits}

\def\GL{\mathop{\rm GL}\nolimits}

\def\Id{\mathop{\rm Id}\nolimits}

\def\Ad{\mathop{\rm Ad}\nolimits}
\def\ad{\mathop{\rm ad}\nolimits}
\def\det{\mathop{\rm det}\nolimits}

\def\log{\mathop{\rm log}\nolimits}
\def\ev{\mathop{\rm ev}\nolimits}
\def\odd{\mathop{\rm odd}\nolimits}

\def\SL{\mathop{\rm SL}\nolimits}

\def\sign{\mathop{\rm sign}\nolimits}
\def\cotor{\mathop{\rm Cotor}\nolimits}
\def\Hom{\mathop{\rm Hom}\nolimits}
\def\Tot{\mathop{\rm Tot}\nolimits}

\def\Qb{{\mathbb Q}}
\def\Cb{{\mathbb C}}

\def\Rb{{\mathbb R}}

\def\Zb{{\mathbb Z}}

\def\Ac{{\cal A}}

\def\Fc{{\cal F}}

\def\Hc{{\cal H}}

\def\Lc{{\cal L}}
\def\Mc{{\cal M}}

\def\Sc{{\cal S}}
\def\Uc{{\cal U}}

\def\Kc{{\cal K}}
\def\Lc{{\cal L}}
\def\Cc{{\cal C}}

\def\rc{{\cal r}}

\def\a{\alpha}
\def\b{\beta}
\def\d{\delta}
\def\D{\Delta}
\def\g{\gamma}
\def\G{\Gamma}

\def\s{\sigma}

\def\vp{\varphi}
\def\nr{\natural}

\def\0b{\bf 0}

 \def\bbc{\mathbb{C}}
\def\bash{\backslash}

\def\fl{\forall}
\def\ify{\infty}

\def\ot{\otimes}

\def\ra{\rightarrow}

\def\sbs{\subset}

\def\rt{\rtimes}
\def\lt{\triangleleft}

\def\bi{\blacktriangleright\hspace{-4pt}\vartriangleleft }
\def\hd{\overset{\ra}{\partial}}
\def \vd{\uparrow\hspace{-4pt}\partial}
\def\hs{\overset{\ra}{\sigma}}
\def \vs{\uparrow\hspace{-4pt}\sigma}
\def\hta{\overset{\ra}{\tau}}
\def \vta{\uparrow\hspace{-4pt}\tau}
\def\hb{\overset{\ra}{b}}
\def \vb{\uparrow\hspace{-4pt}b}
\def \dvb{\uparrow b}
\def\hB{\overset{\ra}{B}}
\def \vB{\uparrow\hspace{-4pt}B}
\def\p{\partial}

\def\0D{\Delta^{(0)}}
\def\1D{\Delta^{(1)}}

\def\wt{\widetilde}
\def\td{\tilde}

\newcommand{\FC}{\mathfrak{C}}
\newcommand{\FD}{\mathfrak{D}}

\newcommand{\Fa}{\mathfrak{a}}
\newcommand{\Fg}{\mathfrak{g}}

\newcommand{\FX}{\mathfrak{X}}
\newcommand{\FY}{\mathfrak{Y}}

\newtheorem{theorem}{Theorem}[section]

\newtheorem{proposition}[theorem]{Proposition}
\newtheorem{lemma}[theorem]{Lemma}
\newtheorem{corollary}[theorem]{Corollary}

\newtheorem{definition}[theorem]{Definition}

\def\ni{\noindent}

\def\build#1_#2^#3{\mathrel{
\mathop{\kern 0pt#1}\limits_{#2}^{#3}}}
\newcommand{\ps}[1]{~\hspace{-4pt}_{^{(#1)}}}
\newcommand{\ns}[1]{~\hspace{-4pt}_{_{{<#1>}}}}
\newcommand{\sns}[1]{~\hspace{-4pt}_{_{{<\overline{#1}>}}}}

\def\odots{\ot\dots\ot}
\newcommand{\nm}[1]{{\mid}#1{\mid}}

\numberwithin{equation}{section}

\begin{document}

\title{\bf Cyclic cohomology of Hopf algebras of transverse symmetries in
codimension $1$}

\author{Henri Moscovici\thanks{Research
    supported by the National Science Foundation
    award no. DMS-0245481.}  \quad and \quad
      Bahram Rangipour   \\ \\
    Department of Mathematics \\
    The Ohio State University \\
    Columbus, OH 43210, USA
}

\date{ \ }

\maketitle

\begin{abstract} 
We develop intrinsic tools for computing  the periodic Hopf cyclic cohomology of Hopf algebras
related to transverse symmetry in codimension 1. Besides the Hopf algebra
found by Connes and the first author in their work on the local index formula
for transversely hypoelliptic  operators on foliations,
this family includes its `Schwarzian' quotient, on which
the Rankin-Cohen universal deformation formula is based, the extended
Connes-Kreimer Hopf algebra related to
renormalization of divergences in QFT, as well as a
series of cyclic coverings of these
Hopf algebras, motivated by the treatment of transverse
symmetry for nonorientable foliations.

The method for calculating their Hopf cyclic cohomology
 is based on two computational devices, which work in tandem and complement
each other: one is a spectral sequence for bicrossed product
Hopf algebras and the other a Cartan-type homotopy formula in Hopf cyclic cohomology.
 \end{abstract}

\section*{Introduction}

Connes' general definition of cyclic cohomology as a derived functor over the
cyclic category~\cite{Cext} allows the flexibility of applying it well beyond
the realm of algebras, for which it was originally formulated~\cite{Cobw, cng} as
a noncommutative de Rham theory.
In particular, the process of computing the local index formula
for transversely hypoelliptic  operators on foliations has led to a
cyclic cohomology apparatus specific to Hopf algebras~\cite{cm2}.
Unfolding the algebraic meaning of the iterated
commutator operations involved in the local index formula~\cite{cm1}
revealed the existence
of a `universal' Hopf algebra $\Hc_n$, that plays for the normal bundle
to a codimension $n$ foliation a role analogous to that of
the affine group of the frame bundle
to a manifold of dimension $n$. Moreover, the calculation itself turned out
to be governed by the cyclic cohomology of a specific type of cyclic module
structure associated to Hopf algebras. The solution provided in~\cite{cm2}
to the transverse index problem was then based on proving
the existence of an explicit isomorphism between the
periodic cyclic cohomology of the cyclic module associated
to the Hopf algebra $\Hc_n$
and the Gelfand-Fuks cohomology of the Lie
algebra ${\Fa}_n$ of formal vector fields on $\Rb^n$, thus
making the link with the `classical' theory of characteristic classes of
foliations.
\medskip

The proof given in~\cite{cm2} to the above
isomorphism was constructive and shed valuable
light  on how to perform explicit
computations in the realm of Hopf cyclic cohomology, but from a technical
viewpoint the
end result was still imported from Lie algebra cohomology.
It thus remained highly desirable to develop intrinsic methods, apt to directly
compute the cyclic cohomology of Hopf algebras such as $\Hc_n$. Furthermore,
one could reasonably expect these methods to require considerably
 less computational effort in the codimension $1$ case.

\medskip

In this paper we present a direct
method for computing
the periodic Hopf cyclic cohomology of a whole class of
Hopf algebras
related to transverse symmetry in codimension $1$.
Besides  $\Hc_1$, this class contains a series of Hopf algebras whose periodic Hopf
cyclic cohomology was not previously calculated.
This family includes
the `Schwarzian' quotient $\Hc_{1 s}$ of $\Hc_1$, which plays a pivotal role in
the universal deformation formula~\cite{cm5'} based on the Rankin-Cohen 
brackets on modular forms, as well as
the extended version of the
Connes-Kreimer Hopf algebra~\cite{ck}
related to the renormalization of divergences in QFT. It also contains
a series of
`cyclic coverings' of these Hopf algebras, obtained as coextensions of them
 by group rings of cyclic groups.
\medskip 
 
The cyclic coverings of $\Hc_1$ are introduced in $\S 1$,
starting from the observation that the process of extension of
the Hopf algebra symmetries to the case of transversely
nonorientable foliations leads naturally to the construction of
a certain `double cover' $\Hc_1^{\dag|2}$ (see $\S 1.1$). A
similarly defined  `infinite cyclic
cover' $\Hc_1^\dag$  (see $\S 1.2$) is then shown to act effectively
on Hecke algebras related to modular forms, of the type introduced
in~\cite{cm5}  (see $\S 1.3$).
\medskip

The tool kit we employ in order to calculate the periodic cyclic cohomology
of the above family of Hopf algebras consists of two computational devices,
which fortunately, albeit rather fortuitously, complement each other. 
The first is
the construction of a bicocyclic module for a class of bicrossed product
Hopf algebras of the type mentioned above,
that allows the computation of their Hopf cyclic cohomology 
(see Theorems~\ref{mix} and \ref{mix+}) by means of
an Eilenberg-Zilber theorem for byparacyclic modules, cf.~\cite{ ft, gj}.
 The second is a Cartan homotopy
formula for the Hopf cyclic cohomology of coalgebras with coefficients
in SAYD modules (Proposition~\ref{homotop}).  
Neither one of these tools would suffice by itself to achieve the desired
goal. However, when applied in tandem they allow the calculation
of the periodic Hopf cyclic cohomology of
$\Hc_1$ (Theorem~\ref{HP(H_1)}) and
of the Schwarzian factor Hopf algebra $\Hc_{1 s}$
(Theorem~\ref{HP(H_s)}) without relying on
Gelfand-Fuks cohomology. Moreover,  they also allow us to compute
the periodic Hopf cyclic cohomology of the cyclic cover
$\Hc_1^\dag$ (with coefficients in various
 modular pairs in involution, cf.~Theorem~\ref{final1}),
as well as to calculate
the periodic Hopf cyclic cohomology of the extended Connes-Kreimer
Hopf algebra $\Hc_{CK}$ (Theorem~\ref{final2})
and of its own cyclic covers (Theorem~\ref{final3}).

\tableofcontents

\section{Coverings of $\Hc_1$ \label{S1}
by cyclic groups} \label{secext}

After briefly reviewing some basic facts about the Hopf algebra $\Hc_1$
(cf. \cite{cm2}),
we shall explain how the consideration of the transversely nonorientable
case leads naturally to the introduction of a `double cover' $\Hc_1^{\dag|2}$
and, more generally, of an infinite cyclic cover $\Hc_1^\dag$. The latter
will be shown to act
 effectively as `symmetries' of modular Hecke algebras (comp. \cite{cm5}).
\medskip

\subsection{$\Hc_1$ and its periodic classes} \label{S1.1}
Given a discrete subgroup of orientation
preserving diffeomorphisms $\G \sbs \Diff^+ (S^1)$, let
 $ F^{+} S^1$ denote the bundle of positively oriented frames
over $S^1$,  on which the  diffeomorphisms $\vp \in \G$
act in the obvious natural manner:
 \begin{equation*}
{\wt \vp} (x, y) = (\vp (x) , \, \vp' (x) \cdot y) \, , \qquad
(x, y) \in F^+ S^1 \simeq (\Rb/\Zb) \times {\Rb}^+ \, .
\end{equation*}
One then forms the crossed product algebra
 \begin{equation*}
 \Ac^+_{\G} \, = \, C_c^{\ify} (F^+ S^1) \rtimes \G \, ,
\end{equation*}
which can be viewed as the linear span of the set of monomials
$\{  f U_{\vp} \, ; \,
 f \in C_c^{\ify} (F^+ S^1) ,  \, \vp \in \G \}$,
 endowed with the product
\begin{equation*}
  f U_{\vp} \cdot g U_{\psi} = f \, (g \circ \vp^{-1}) \, U_{\vp \psi} \, .
\end{equation*}
The algebra $\Ac^+_{\G}$ typifies the `coordinates' of the `space of leaves'
of a \textit{transversely orientable} codimension $1$ foliation.
\medskip

 The trivial connection on the frame bundle $\,  F^+S^1 \ra S^1$ is
 implemented by the vector fields
\begin{equation}  \label{conn}
Y \, = \,  y \, \frac{\partial }{\partial y} \qquad \text{and}
  \qquad X \, = \,  y \, \frac{\partial}{\partial x} \, ,
 \end{equation}
 the first being the generator of the vertical sub-bundle and the second
 the generator of the horizontal sub-bundle.
These two basic vector fields become linear operators on $\Ac^+_{\G}$,
acting as
\begin{equation}  \label{yx}
Y (f \, U_{\vp}) = Y ( f) \, U_{\vp} \, , \qquad
X (f \, U_{\vp}) = X ( f) \, U_{\vp} \, .
\end{equation}
The prolongation of $Y$ acts as derivation
\begin{equation} \label{hy}
Y (ab) = Y (a) \, b + a \, Y (b) \, , \qquad a,b \in \Ac^+_{\G} \, ,
\end{equation}
because $Y$ is invariant under the action of $ \Diff^+ (S^1)$ on
$ F^+S^1$. On the other hand,
\begin{equation}   \label{xphi}
X (f \circ \vp) (x, y) =X ( f) (\vp (x) , \vp' (x) \cdot y)
+ y  \frac{\vp^{''}(x)}{\vp^{'}(x)} Y (f)  (\vp (x) ,  \vp' (x) \cdot y)
\end{equation}
which is tantamount with the following `Leibniz rule' for its prolongation:
\begin{equation} \label{hx}
X (ab) = X (a) \, b \, + \, a \, X (b) + \d_1 (a) \, Y (b)  \, ,
\end{equation}
where
\begin{equation} \label{d1}
 \d_1 (f \, U_{\vp^{-1}}) \, = \, y \,  \frac{d}{d x}
\left(\log \frac{d\vp}{dx}\right) \, f U_{\vp^{-1}}  \, .
\end{equation}
In turn, the operator  $\d_1$ is itself a derivation,
\begin{equation} \label{hd1}
 \d_1 (ab) = \d_1 (a) \, b + a \, \d_1 (b) \, ,
\end{equation}
and its successive
commutators with $X$ produces new operators
\begin{equation} \label{dn}
 \d_n (f \, U_{\vp^{-1}}) \, = \,  y^n \,  \frac{d^{n}}{d x^{n}}
\left(\log \frac{d\vp}{dx}\right) \, f U_{\vp^{-1}} \, , \qquad \fl \, n \geq 1 \, ,
\end{equation}
that satisfy progressively more complicated Leibniz rules.
\smallskip

The formulae \eqref{hy}, \eqref{hx}, \eqref{hd1} in conjunction
with \eqref{yx},  \eqref{d1}, \eqref{dn} express precisely the
fact that they define a Hopf action of the Hopf algebra
 $\Hc_1$ on $\Ac^+_{\G}$. As an algebra
$\Hc_1$ is
generated by $X$, $Y$, $\delta_k$,  $k\in \mathbb{N}$,
subject to  the relations
\begin{equation} \label{pres}
[Y, X]=X,\qquad
[Y,\delta_k]=k\delta_k,\qquad
[X,\delta_k]=\delta_{k+1}, \qquad
[\delta_j, \delta_k]=0.
\end{equation}
Its coalgebra structure is uniquely determined by
\begin{align} \label{copr1}
&\Delta(Y)=Y\ot 1+1\ot Y,\\ \label{copr2}
&\Delta(X)=X\ot 1+1\ot X+\delta_1\ot Y,\\ \label{copr3}
& \Delta(\delta_1)=\delta_1\ot 1+1\ot \delta_1,\\ \nonumber \\ \label{coun}
&\epsilon(X)=\epsilon(Y)=\epsilon(\delta_k)=0,
\end{align}
together with the requirement that $\Delta : \Hc_1 \ra \Hc_1 \ot \Hc_1$
is multiplicative.

The antipode is the unique antihomomorphism $S : \Hc_1 \ra \Hc_1$
defined on generators by
\begin{align} \label{antip}
S(X)=-X+\delta_1Y, ~~~S(Y)=-Y ~~~\text{ and}~~~
S(\delta_1)=-\delta_1\, .
\end{align}
\medskip

Furthermore, the above action admits an invariant trace.
Indeed, the volume form $\displaystyle \frac{dx\wedge
dy}{y^2} \, $ on $ \, F^+ S^1 \,$ is invariant under $\,
\Diff^+ (S^1) \,$ and thus gives rise to a trace $\, \tau^+ :
\Ac^+_{\G} \ra \mathbb{C} $,
\begin{equation} \label{tr+}
\tau^+ (f U_{\vp}) \, = \, \begin{cases} \displaystyle
 \int_{ F^+ S^1} f(x ,y) \, \frac{dx\wedge dy}{y^2} \, ,  \qquad
 &\text{ if }  \vp = 1 \, ,  \\ \\
\qquad 0 \, , \qquad &\text{ if }  \vp \ne 1 \, .
\end{cases}
\end{equation}
This trace satisfies the invariance property
\begin{equation} \label{invt}
\tau^+ (h(a)) = \d (h) \, \tau^+ (a) \, , \qquad \fl \quad h \in
\Hc_{1} \, ,
\end{equation}
where $\,\d \in \Hc_1^* \,$ is  the `modular' character, defined
on generators by
\begin{equation} \label{mch}
\d (Y) = 1, \quad \d (X) = 0 , \quad \d (\d_n) = 0 \, .
\end{equation}
To prove the identity \eqref{invt}, it suffices to note that it
holds for  $X$ and $Y$ by an elementary integration by parts
argument, and is trivially fulfilled by $\d_1$ since $\d_1 (\Id) = 0$.
\bigskip

For a quick exposition of the basic notions and notation regarding Hopf cyclic
cohomology,  the reader is referred to \cite{cm3} (see also $\S$\ref{S2.1} below).
In particular, a {\em modular pair}  for a Hopf algebra
$\Hc$ is a pair  $(\sigma, \delta)$, where $\sigma$ is a group-like
element of $\Hc$ and  $\, \delta:\Hc \rightarrow \bbc\, $ is an algebra map
such that $\d(\sigma)=1$. It is  called {\em modular pair in
involution}, to be abbreviated henceforth {\em MPI},
 iff denoting $\, S_\delta = \delta\ast S$, the following condition is
 satisfied:
 \begin{equation} \label{mpi}
   S_\delta^2 (h)=\sigma h\sigma^{-1}\, , \quad \fl \, h \in \Hc\, .
 \end{equation}

The pair
($\d,1$) forms such an MPI,  so that
the corresponding Hopf cyclic cohomology $HC_{(\d, 1)}^* (\Hc_1)$  is well-defined.
Moreover, the assignment
\begin{equation*}
\chi^n_{\tau^+} (h^{1} \ot \ldots \ot h^{n}) (a^0 , \ldots , a^{n})
=  \tau^+ (a^0 \, h^{1} (a^{1}) \ldots h^{n} (a^{n})) \, ,
\quad h^i \in \Hc_1 , \, \, a^j \in \Ac^+_{\G}
\end{equation*}
 induces a characteristic homomorphism in cyclic cohomology (cf. \cite{cm2, cm3})
\begin{equation} \label{charm}
 \chi_{\tau^+}^{*} \, : \,  HC_{(\d, 1)}^*  (\Hc_1) \, \ra \, HC^* \,
(\Ac^+_{\G}) \, .
\end{equation}
\medskip

 As is well-known, in codimension $1$ there are just two linearly independent
 geometric characteristic classes of foliations: the $0$-dimensional Pontryagin
 class $p_0 = 1$ and the Godbillon-Vey class. In the framework of cyclic cohomology,
 the first becomes the \textit{transverse fundamental class} and
 is represented by the cyclic $2$-cocycle $TF_\G \in  ZC^2 (\Ac^+_{\G}) $,
 \begin{equation} \label{tfc}
TF_\G (f^0 U_{\vp_0^{-1}} , f^1 U_{\vp_1^{-1}} , f^2 U_{\vp_2^{-1}}) =
\left\{ \begin{matrix}  \displaystyle
\int_{F^+S^1} f^0  {\wt \vp_0}^* ( d f^1)  {\wt \vp_0}^* {\wt \vp_1}^* (d f^2),
\, \, \,
 \vp_2  \vp_1 \vp_0 = \Id \cr\cr
\qquad 0 \, , \qquad \qquad \vp_2  \vp_1 \vp_0 \neq \Id \, ,
\end{matrix} \right.
\end{equation}
and the second by the cyclic $1$-cocycle
 \begin{equation} \label{gv}
GV_\G (f^0 U_{\vp_0^{-1}} , f^1 U_{\vp_1^{-1}} ) =
\left\{ \begin{matrix}  \displaystyle
\int_{F^+S^1}  f^0 \, {\wt \vp_0}^* (f^1) \, y \frac{\vp_0''(x)}{\vp_0'(x)} \,
\frac{dx \wedge dy}{y^{2}}, \quad
 \vp_1  \vp_0 = \Id \cr\cr
\qquad 0 \, ,  \qquad \qquad \vp_1 \vp_0 \neq \Id \, ,
\end{matrix} \right.
\end{equation}
\smallskip

By an elementary calculation (see \cite[Prop. 18]{cm6}), the first cocycle
can be expressed in the form
\begin{eqnarray*}
&\,& TF_\G (f^0 U_{\vp_0^{-1}} , f^1 U_{\vp_1^{-1}} , f^2 U_{\vp_2^{-1}})\, =\,
 \tau^+ (f^0 U_{\vp_0^{-1}} \, X(f^1 U_{\vp_1^{-1}} )\, Y( f^2 U_{\vp_2^{-1}}) \\
 &-&\tau^+ (f^0 U_{\vp_0^{-1}} \, Y(f^1 U_{\vp_1^{-1}} )\, X( f^2 U_{\vp_2^{-1}})
- \tau^+ (f^0 U_{\vp_0^{-1}} \, \d_1Y(f^1 U_{\vp_1^{-1}} )\, Y( f^2 U_{\vp_2^{-1}})
\end{eqnarray*}
which can be recognized as  the image
by the characteristic map \eqref{charm}
of the Hopf cyclic cocycle $TF \in ZC^2 _{(\d, 1)}  (\Hc_1) $
\begin{equation} \label{hotfc}
 TF \,= \, X \ot Y - Y \ot X - \d_1 \,  Y \ot Y \, \in \, \Hc_1 \ot \Hc_1\, .
\end{equation}
\smallskip

Using the invariance property \eqref{invt} for $h=\d_1$,
the Godbillon-Vey cocycle \eqref{gv} can be rewritten as
\begin{eqnarray*}
GV_\G (f^0 U_{\vp_0^{-1}} , f^1 U_{\vp_1^{-1}} )&=&
\tau (\d_1 (f^0 U_{\vp_0^{-1}}) \,  f^1 U_{\vp_1^{-1}} )
=- \tau (  f^0 U_{\vp_0^{-1}}  \, \d_1 (f^1 U_{\vp_1^{-1}}) ) \\
&=& - \chi^1_{\tau} (\d_1) (f^0 U_{\vp_0^{-1}} ,  f^1 U_{\vp_1^{-1}}) \, ,
\end{eqnarray*}
which is obviously the image
by the characteristic map of the Hopf cyclic cocycle
 \begin{equation*}
 GV \,= \, - \d_1 \, \in \, ZC^1 _{(\d, 1)} (\Hc_1)   \, .
\end{equation*}

\subsection{$\Hc_1^\dag$ and its periodic classes} \label{S1.2}

We now proceed to upgrade the above constructions in order to cover the case of
non-orientable codimension $1$ foliations. First of all, this entails replacing the
positive frame bundle $F^+S^1$ by the quotient $FS^1/ \Zb_2$ of the
full frame bundle, where the action of $\Zb_2$ is given by the reflection
$(x, y) \mapsto (x, -y)$,
and allowing $\G$ to be an arbitrary subgroup of the full diffeomorphism
group $\Diff (S^1)$. The corresponding crossed product algebra is then
 \begin{equation*}
 \Ac_{\G} \, = \, C_c^{\ify} (F S^1/ \Zb_2) \rtimes \G \, ,
\end{equation*}
where $ \, C_c^{\ify} (F S^1/ \Zb_2)$ is identified with the space of even
functions
\begin{equation*}
C_c^{\ify} (F S^1)^{\ev} \, = \, \{ f \in  C_c^{\ify} (F S^1) \, \, ; \quad f (x, y) =  f (x, -y) \, ,
\quad \fl \,  (x, y) \in S^1 \times \Rb \} \, .
\end{equation*}
The trivial connection is now implemented by the vector fields $\{Y, X\}$,
where $Y$ is the same as in \eqref{conn} while the expression
 of $X$ needs an obvious adjustment:
\begin{equation}  \label{nconn}
 Y \, = \,  y \, \frac{\partial}{\partial y} \qquad \text{and} \qquad
 X \, = \,  |y| \, \frac{\partial}{\partial x} \, .
 \end{equation}
As a consequence, \eqref{xphi} becomes
\begin{eqnarray}  \nonumber
X (f \circ \vp) (x, y) &=& \sign \vp' (x)\cdot X ( f) (\vp (x) , \vp' (x) \cdot y) \\  \label{oxphi}
&+& |y|  \frac{\vp^{''}(x)}{\vp^{'}(x)} Y (f)  (\vp (x) ,  \vp' (x) \cdot y) \, ,
\end{eqnarray}
which implies the Leibniz rule
\begin{equation} \label{ohx}
X (ab) = X (a) \, b \, + \, \s (a )\, X (b) + \d_1 (a) \, Y (b)  \, ,
\end{equation}
where
\begin{eqnarray} \label{os}
\s (f \, U_{\vp^{-1}}) &=& \sign \vp' (x) \cdot f \, U_{\vp^{-1}}  \, , \\  \label{od1}
\text{and} \qquad  \d_1 (f \, U_{\vp^{-1}}) &=& |y| \, \frac{d}{d x}
\left(\log \frac{d\vp}{dx}\right) \, f U_{\vp^{-1}}  \, .
\end{eqnarray}
The operator  $\d_1$ also becomes an $\s$-derivation,
\begin{equation} \label{ohd1}
 \d_1 (ab) = \d_1 (a) \, b + \, \s (a )\, \d_1 (b) \, ,
\end{equation}
while $\s$ acts as automorphism
\begin{equation} \label{as}
 \s (ab) \, = \, \s(a) \, \s(b)  \, .
\end{equation}
The iterated commutators of $\d_1$ with
 $X$ give rise to the operators
\begin{equation} \label{odn}
 \d_n (f \, U_{\vp^{-1}}) \, = \,  |y|^n \,  \frac{d^{n}}{d x^{n}}
\left(\log \frac{d\vp}{dx}\right) \, f U_{\vp^{-1}} \, , \qquad \fl \, n \geq 1 \, .
\end{equation}
\smallskip

The formulae \eqref{hy}, \eqref{ohx},  \eqref{ohd1} and \eqref{as} lead to
the definition of a `central extension'
of $\Hc_1$ by the group ring
$\, \Cb [\Zb] \equiv \Cb [ \s, \s^{-1}]$,
\begin{equation} \label{Hdag}
\Hc_1^\dag \, : = \, \Hc_1 \ot \Cb [ \s, \s^{-1}] \, ,
\end{equation}
whose presentation as an algebra is given
by \eqref{pres} plus the condition that $\s$ is central. The coalgebra structure
of $\Hc_1^\dag$ is determined by obvious analogues of the
formulae \eqref{copr1}--\eqref{coun}, namely
\begin{align*}
&\Delta(Y)=Y\ot 1+1\ot Y,\\
&\Delta(X)=X\ot 1+ \s \ot X+\delta_1\ot Y,\\
& \Delta(\delta_1)=\delta_1\ot 1+\s \ot \delta_1,\\
& \Delta(\s)=\s \ot \s , \quad
\Delta(\s^{-1})=\s^{-1} \ot \s^{-1} \\ \nonumber \\
&\epsilon(X)=\epsilon(Y)=\epsilon(\delta_k)=0, \quad \epsilon (\s) =1,
\end{align*}
while the corresponding antipode is defined by the relations
\begin{eqnarray*}
S(\s)&=&\s^{-1} \, , ~~~~
S(X)=\s^{-1}(-X+\delta_1Y) \, , \\
S(Y)&=&-Y\, , ~~~~
S(\delta_1)=-\s^{-1}\delta_1\, .
\end{eqnarray*}
\medskip

For each
integer $N \geq 1$ one can add to the presentation of  $\Hc_1^\dag$ the
relation
\begin{equation*}
 \s^N  \, = \, 1  \, ,
\end{equation*}
or equivalently replace $ \Cb [\Zb] $ by $ \Cb [\Zb/N \Zb] $, and preserve
the rest of the axioms  to define
a Hopf algebra $\Hc_1^{\dag|N}$. Alternatively,
$\Hc_1^{\dag|N}$ can be defined as the quotient of $\Hc_1^{\dag}$
by the ideal generated by $\, \s^N  - 1$, which is also a coideal.
\medskip

At this point we note that, with the notation just introduced,
 the identities \eqref{yx}, \eqref{os},  \eqref{od1} and \eqref{odn}
actually define a Hopf action of
 $\Hc_1^{\dag|2}$ on $\Ac_{\G}$. Furthermore, the obvious counterpart of
the trace $\tau^+$, namely
\begin{equation*}
\tau (f U_{\vp}) \, = \, \begin{cases} \displaystyle
 \int_{ F S^1/ \Zb_2} f(x ,y) \, \frac{dx \wedge dy}{y|y|} \, ,  \qquad
 &\text{ if }  \vp = 1 \, ,  \\ \\
\qquad 0 \, , \qquad &\text{ if }  \vp \ne 1 \, ,
\end{cases}
\end{equation*}
gives an invariant $\s^{-1}$-trace on $\Ac_{\G}$, which means that it satisfies
the identity
\begin{equation*}
\tau(a\, b) \, = \, \tau (b\, \s^{-1} (a)) \, , \qquad \fl \, a, b \in \Ac_{\G}
\end{equation*}
Indeed, the change of variables $(x,|y|) \mapsto (\vp(x), |\vp'(x)| |y|)$ gives
\begin{eqnarray*}
\tau (f U_{\vp} \cdot g U_{\vp^{-1}}) &=&\int_{F S^1/ \Zb_2} f(x ,y) \, (g \circ \tilde\vp^{-1} (x, y))
\, \frac{dx \wedge dy}{y|y|}  \\ \\
&=& \int_{F S^1/ \Zb_2} g(x ,y) \, (f \circ \tilde\vp (x, y)) \sign \vp' (x)
\, \frac{dx \wedge dy}{y|y|}  \\ \\
&=&\tau (\s(g U_{\vp^{-1}}) \, f U_{\vp} )\, = \, \tau (g U_{\vp^{-1}} \cdot \s^{-1} (f U_{\vp}) ) \, ;
\end{eqnarray*}
in the last line we used the invariance property
\begin{equation*}  \label{invs}
 \tau(\s(a)) \, = \,   \tau (a) \, , \qquad \fl \quad a \in
\Ac_{\G} \, ,
\end{equation*}
which is an immediate consequence of the
fact that $\s (\Id) = 1$.
\smallskip

More generally, as in the case of  $\tau^+$ one can easily
check that
\begin{equation}  \label{invst}
 \tau(h(a)) \, = \,  \d (h) \, \tau (a) \, , \qquad \fl \quad h \in
\Hc_{1} \, ;
\end{equation}
here  the modular character $\d$ is the extension of that defined
by \eqref{mch}, with the additional condition
\begin{equation} \label{omch}
  \d (\s) = 1 \, .
\end{equation}
As $(\d, \s^{-1})$ is obviously an involutive modular pair for $\Hc_1^{\dag}$,
one obtains a characteristic homomorphism in cyclic cohomology
\begin{equation} \label{ocharm}
 \chi_{\tau}^{*} \, : \, HC_{(\d, \s^{-1})}^* (\Hc_1^{\dag}) \, \ra HC^*(\Ac_\G)
\end{equation}
that factors through $\Hc_1^{\dag|2}$,
by means of the assignment
\begin{equation*}
\chi^n_{\tau} (h^{1} \ot \ldots \ot h^{n}) (a^0 , \ldots , a^{n})
=  \tau (a^0 \, h^{1} (a^{1}) \ldots h^{n} (a^{n})) ,
\quad h^i \in \Hc_1^{\dag} , \, \, a^j \in \Ac_\G .
\end{equation*}
\bigskip

From the very definitions of the two characteristic maps,
it follows quite easily that if
$\G \sbs \Diff (S^1)$ is a discrete subgroup of diffeomorphisms,
 $\G^+ = \G \cap \Diff^+ (S^1)$, and
  $ \iota_{\G}  : \Ac^+_{\G^+} \ra \Ac_\G$ denotes
  the obvious inclusion, then the diagram
  \begin{equation*}
 \begin{xy}
\xymatrix{  HC_{(\d, \s^{-1})}^* (\Hc_1^{\dag})   \ar[r]^{\pi_1^*}
\ar[d]^{ \chi_{\tau}^{*}}&  HC_{(\d, 1)}^* (\Hc_1)  \ar[d]^{ \chi_{\tau^+}^{*}} \\
HC^* (\Ac_\G) \ar[r]^{\iota_{\G}^*} &  {HC^* (\Ac^+_{\G^+} ) }  }
\end{xy}
\end{equation*}
is commutative; here
$ \pi_1 : \Hc_1^{\dag} \ra   \Hc_1$ stands for the natural projection
 that sends $\s$ to $1$ and $ {\pi_1}^*$ for the map induced on cohomology.
\medskip

The two classes given by the cyclic cocycles \eqref{tfc} and \eqref{gv}
 lift naturally to classes in $\, HC^* (\Ac_\G)$:
 \begin{equation} \label{otfc}
TF_\G (f^0 U_{\vp_0^{-1}} , f^1 U_{\vp_1^{-1}} , f^2 U_{\vp_2^{-1}}) =
\left\{ \begin{matrix}  \displaystyle
\int_{FS^1/\Zb_2} f^0  {\wt \vp_0}^* ( d f^1)  {\wt \vp_0}^* {\wt \vp_1}^* (d f^2),
\, \,
 \vp_2  \vp_1 \vp_0 = \Id \cr\cr
\qquad 0 \, , \qquad \qquad \vp_2  \vp_1 \vp_0 \neq \Id \, ,
\end{matrix} \right.
\end{equation}
respectively
 \begin{equation} \label{ogv}
GV_\G (f^0 U_{\vp_0^{-1}} , f^1 U_{\vp_1^{-1}} ) =
\left\{ \begin{matrix}  \displaystyle
\int_{FS^1/\Zb_2}  f^0 \, {\wt \vp_0}^* (f^1) \, |y| \frac{\vp_0''(x)}{\vp_0'(x)} \,
\frac{dx \wedge dy}{y|y|} , \quad
 \vp_1  \vp_0 = \Id \cr\cr
\qquad 0 \, ,  \qquad \qquad \vp_1 \vp_0 \neq \Id \, ,
\end{matrix} \right.
\end{equation}

 As in the orientable case, we can identify them as
 characteristic images of classes in
 $HC_{(\d, \s^{-1})}^* (\Hc_1^{\dag})  $.
 \medskip

 \begin{proposition}
 The elements
 \begin{eqnarray*}
 GV^\dag&=& - \s^{-1} \d_1   \, , \\
 TF^\dag&=&\s^{-1} X \ot \s^{-1} Y - Y \ot \s^{-1} X
- \s^{-1}\d_1 Y \ot \s^{-1} Y
 \end{eqnarray*}
 are cyclic cocycles
 in the Hopf cyclic module associated to the involutive modular pair
 $(\Hc^{\dagger}_1; \d , \s^{-1})$, and they
 satisfy the naturality property
 \begin{eqnarray*}
 \chi_\tau^* (GV^\dag)&=&GV_\Gamma \, , \qquad {\pi_1^2}^* (GV^\dag)= GV \, ,
 \\
 \chi_\tau^* (TF^\dag) &=&TF_\Gamma \, , \qquad {\pi_1^2}^* (TF^\dag)= TF \, .
 \end{eqnarray*}
    \end{proposition}
 \medskip

 \begin{proof}
 Starting with the Godbillon-Vey cocycle, we note that formula \eqref{ogv}
 is equivalent to
 \begin{equation*}
 GV_\G (f^0 U_{\vp_0^{-1}} , f^1 U_{\vp_1^{-1}} )\, = \,
\tau (\d_1 (f^0 U_{\vp_0^{-1}}) \,  f^1 U_{\vp_1^{-1}} ) \, ,
 \end{equation*}
 and using the invariance property \eqref{invst}, for $\d_1$ and then for
$\s$, it follows that
 \begin{eqnarray*}
&\,&GV_\G (f^0 U_{\vp_0^{-1}} , f^1 U_{\vp_1^{-1}} )\, = \,
 - \tau ( \s( f^0 U_{\vp_0^{-1}})  \,  \d_1 (f^1 U_{\vp_1^{-1}}) ) =\\
&\,& \quad = - \tau (f^0 U_{\vp_0^{-1}}  \, \s^{-1}\d_1 (f^0 U_{\vp_0^{-1}}))
\, = \, - \chi_{\tau} (\s^{-1} \d_1)(f^0 U_{\vp_0^{-1}} , f^1 U_{\vp_1^{-1}} )\, ,
 \end{eqnarray*}
which expresses it as the image
by the characteristic map \eqref{ocharm}
of
 \begin{equation*}
 GV^\dag \,= \, - \s^{-1} \d_1 \, \in \, C^1_{(\d, 1)} ( \Hc_1)   \, .
\end{equation*}
The latter is a Hochschild cocycle, since
  \begin{equation*}
 b (\s^{-1} \d_1) = 1 \ot \s^{-1} \d_1 - (\s^{-1} \ot \s^{-1} )
 (\d_1 \ot 1 + \s \ot \d_1)
 + \s^{-1} \d_1 \ot \s^{-1} = 0 \, ,
\end{equation*}
 and is also  cyclic, because
 \begin{equation*}
 \tau_1 (-\s^{-1} \d_1) \, = \, - \tilde S (\s^{-1} \d_1)) \,\s^{-1} = \,
   - {\tilde S} (\d_1) \, =  \,  \s^{-1} \d_1 \, .
 \end{equation*}
  \medskip

 Passing to the transverse fundamental class, we note first that
relative to the framing \eqref{nconn} of the cotangent bundle $T^*(FS^1/\Zb_2)$,
for any $f \in C^\infty (FS^1/\Zb_2)$,
\begin{equation*}
df \, = \, X(f) \, |y|^{-1} dx \, + \, Y(f)\,  y^{-1} dy  \, ,
\end{equation*}
hence for any $\vp \in \Diff (S^1)$
\begin{equation*}
 {\wt \vp}^* ( d f) \, =\, d(f \circ \wt\vp) \, = \,
 X(f \circ \wt\vp)\, |y|^{-1} dx \, + \,
  Y(f \circ \wt\vp)\,  y^{-1} dy   \, .
\end{equation*}
From \eqref{oxphi} it follows that
\begin{eqnarray*}
 {\wt \vp}^* ( d f)&=&\left(\sign \vp' (x)\, X ( f)  \circ \wt\vp  \, + \,
  |y|  \frac{\vp^{''}(x)}{\vp^{'}(x)} \, Y (f) \circ \wt\vp\right)
\, |y|^{-1} dx   \\
&+& \, \left(Y(f) \circ \wt\vp \right) \,  y^{-1} dy   \, .
\end{eqnarray*}
On substituting in \eqref{otfc} the above expression applied to
$ {\wt \vp_0}^* ( d f^1)$ as well as  to $ \wt \vp_{01}^*( d f^2)$,
where $ \vp_{01} :=  \vp_{1} \vp_{0}$,  one obtains
 \begin{eqnarray} \label{lhs}
&\,&f^0  {\wt \vp_0}^* ( d f^1)  {\wt \vp_0}^* {\wt \vp_1}^* (d f^2) \, = \\ \nonumber
&=& f^0 \Big( \big(\sign \vp_0' (x)\, X ( f^1)  \circ \wt\vp_0 \, + \,
  |y|  \frac{\vp_0^{''}(x)}{\vp_0^{'}(x)} \, Y (f^1) \circ \wt\vp_0 \big)\,
   Y (f^2) \circ \wt\vp_{01} \\ \nonumber
&-&Y(f^1) \circ \wt\vp_0 \, \big(\sign \vp_{01}' (x)\, X ( f^2)  \circ \wt\vp_{01} \, + \,
  |y|  \frac{\vp_{01}^{''}(x)}{\vp_{01}^{'}(x)} \, Y (f^2) \circ \wt\vp_{01} \big)\Big)
\frac{dx \wedge dy}{ |y|y} \\ \nonumber \\ \label{rhs1}
&=&f^0 \cdot \sign \vp_0' (x)\cdot X ( f^1)  \circ \wt\vp_0 \cdot Y (f^2) \circ \wt\vp_{01}
\cdot \frac{dx \wedge dy}{ |y|y} \, \\\nonumber \\ \label{rhs2}
&-&f^0 \cdot  Y(f^1) \circ \wt\vp_0 \cdot  \sign \vp_{01}' (x)\cdot X ( f^2)  \circ \wt\vp_{01}
\cdot  \frac{dx \wedge dy}{ |y|y} \\ \label{rhs3}\\ \nonumber
&-&f^0 \cdot
\left( |y|  \frac{\vp_{01}^{''}(x)}{\vp_{01}^{'}(x)}\, - \,  |y|  \frac{\vp_0^{''}(x)}{\vp_0^{'}(x)} \right)\cdot
Y(f^1) \circ \wt\vp_0 \cdot Y (f^2) \circ \wt\vp_{1} \wt\vp_{0} \cdot \frac{dx \wedge dy}{ |y|y} \, .
\end{eqnarray}
Assuming $\, \vp_{2}\vp_{1} \vp_{0} = \Id$ and taking the integral, the left hand side
of \eqref{lhs} is just
 \begin{equation*}
 TF_\G (f^0 U_{\vp_0^{-1}} , f^1 U_{\vp_1^{-1}} , f^2 U_{\vp_2^{-1}}) \, =  \,
\int_{FS^1/\Zb_2} f^0  {\wt \vp_0}^* ( d f^1)  {\wt \vp_0}^* {\wt \vp_1}^* (d f^2)  \, .
\end{equation*}
On the other hand, after integration each of the three terms in the right hand side
of the  above can be recognized as being in the image of the characteristic map.
Indeed, the term \eqref{rhs1} is identical to
 \begin{eqnarray*}
&\,&\tau \big(\s (f^0 U_{\vp_0^{-1}} ) \, X (f^1 U_{\vp_1^{-1}}) \, Y (f^2 U_{\vp_2^{-1}})\big)= \\
&\,& \qquad = \,
\tau \big(f^0 U_{\vp_0^{-1}}  \,  \s^{-1} X (f^1 U_{\vp_1^{-1}}) \,  \s^{-1}Y (f^2 U_{\vp_2^{-1}})\big) \\
&\,& \qquad = \,
 \chi_{\tau} (\s^{-1} X  \ot  \s^{-1} Y) (f^0 U_{\vp_0^{-1}} , f^1 U_{\vp_1^{-1}} , f^2 U_{\vp_2^{-1}})  \, .
\end{eqnarray*}
 For the second term,  using
\begin{equation*}
 \sign \vp_{01}' (x) \, =  \, \sign \vp_{1}' (\vp_0(x)) \cdot  \sign \vp_{0}' (x)\, ,
\end{equation*}
one identifies  \eqref{rhs2} as
 \begin{eqnarray*}
&\,&\tau \big(\s (f^0 U_{\vp_0^{-1}} ) \, \s Y (f^1 U_{\vp_1^{-1}}) \, X (f^2 U_{\vp_2^{-1}})\big)= \\
&\,& \qquad = \,
\tau \big(f^0 U_{\vp_0^{-1}}  \, Y (f^1 U_{\vp_1^{-1}}) \, \s^{-1}X (f^2 U_{\vp_2^{-1}})\big)\\
&\,& \qquad = \,
  \chi_{\tau} (Y  \ot  \s^{-1} X) (f^0 U_{\vp_0^{-1}} , f^1 U_{\vp_1^{-1}} , f^2 U_{\vp_2^{-1}})  \, .
\end{eqnarray*}
To treat the last term, we note that
 \begin{eqnarray*}
|y|  \frac{\vp_{01}^{''}(x)}{\vp_{01}^{'}(x)} &-& |y|  \frac{\vp_0^{''}(x)}{\vp_0^{'}(x)}
= |y| \frac{d}{dx}\Big(\log \vp_{01}' (x)\Big) - |y| \frac{d}{dx}\Big(\log \vp_{0}' (x)\Big) \\
&=&|y| \frac{d}{dx}\Big(\log \vp_{1}' (\vp_0(x))\Big)
 = \sign \vp'_0 (x) |y  \vp'_0 (x)| \frac{\vp_{1}^{''}(\vp_0(x))}{\vp_{1}^{'}(\vp_0(x))}
\end{eqnarray*}
hence  \eqref{rhs3} can be expressed as
\begin{eqnarray*}
&\,& -\tau \big( \s (f^0 U_{\vp_0^{-1}} ) \,  \d_1 Y (f^1 U_{\vp_1^{-1}}) \,
 Y (f^2 U_{\vp_2^{-1}})\big)= \\
&\,& \qquad = \,
- \tau \big( f^0 U_{\vp_0^{-1}}  \, \s^{-1} \d_1 Y (f^1 U_{\vp_1^{-1}}) \,
\s^{-1} Y (f^2 U_{\vp_2^{-1}})\big) \, .
\end{eqnarray*}
\medskip

To show that $\,TF^\dag$ is a Hochschild cocycle we compute
$\, b (TF^\dag)$, term by term:
 \begin{eqnarray*}
b (- \s^{-1} X \ot \s^{-1} Y ) &=& - 1\ot \s^{-1} X \ot \s^{-1} Y
 + \s^{-1} X \ot \s^{-1}\ot \s^{-1} Y  \\
 &+& 1 \ot \s^{-1} X \ot \s^{-1} Y
+ \s^{-1} \d_1 \ot \s^{-1} Y \ot \s^{-1} Y \\
 &-& \s^{-1} X \ot \s^{-1} Y \ot \s^{-1} - \s^{-1} X \ot \s^{-1} \ot \s^{-1} Y \\
 &+& \s^{-1} X \ot \s^{-1} Y \ot \s^{-1}
 = \quad \s^{-1} \d_1 \ot \s^{-1} Y \ot \s^{-1} Y \, ,
 \end{eqnarray*}
  \begin{eqnarray*}
 b (Y \ot \s^{-1} X )  &=&  1\ot Y \ot \s^{-1} X
 - Y \ot 1 \ot \s^{-1} X - 1 \ot Y \ot \s^{-1} X \\
 &+& Y \ot \s^{-1} X \ot \s^{-1}
+ Y \ot 1 \ot \s^{-1} X +  Y \ot \s^{-1} \d_1 \ot \s^{-1} Y \\
 &-&  Y \ot \s^{-1} X \ot \s^{-1}  \, = \qquad
Y \ot \s^{-1} \d_1 \ot \s^{-1} Y \, ,
 \end{eqnarray*}
 \begin{eqnarray*}
b ( \s^{-1} \d_1 Y \ot \s^{-1} Y ) &=&  1\ot \s^{-1} \d_1 Y \ot \s^{-1} Y
 - \s^{-1} \d_1 Y \ot \s^{-1}\ot \s^{-1} Y  \\
 &-& \s^{-1} \d_1 \ot \s^{-1} Y \ot \s^{-1} Y
- Y \ot \s^{-1} \d_1 \ot \s^{-1} Y  \\
 &-& 1 \ot \s^{-1} \d_1 Y \ot \s^{-1} Y +
 \s^{-1} \d_1 Y \ot \s^{-1} Y \ot \s^{-1}  \\
 &+& \s^{-1} \d_1 Y \ot \s^{-1} \ot \s^{-1} Y
 - \s^{-1} \d_1 Y \ot \s^{-1} Y \ot \s^{-1} \\
&=& - \s^{-1} \d_1 \ot \s^{-1} Y \ot \s^{-1} Y
- Y \ot \s^{-1} \d_1 \ot \s^{-1} Y ;
\end{eqnarray*}
summing up, one obtains $\qquad b (TF^\dag) \, = \, 0$.
\medskip

Finally, we check the cyclicity of $TF^\dag$.  One has
 \begin{eqnarray*}
 \tau_2 (S(X) \ot \s^{-1}Y) &=& S_\d(S(X))_{(1)} \s^{-1}Y \ot
 S_\d(S(X))_{(2)}) \s^{-1} \\
  &=& X \s^{-1} Y \ot \s^{-1} +
  Y \ot X \s^{-1} +  \d_1 \s^{-1} Y \ot Y \s^{-1}  \, ,
\end{eqnarray*}
\begin{eqnarray*}
 \tau_2 (Y \ot \s^{-1}X) &=& S(Y_{(2)}) \s^{-1} X \ot
 S_\d(Y_{(1)}) \s^{-1} \\
  &=& \s^{-1} X \ot (- Y + 1) \s^{-1} -
  Y \s^{-1} X \ot \s^{-1} \\
   &=& - \s^{-1} X \ot Y \s^{-1} + \s^{-1} X \ot \s^{-1} -
  Y \s^{-1} X \ot \s^{-1} ;
\end{eqnarray*}
adding up and taking into account that
\begin{eqnarray*}
&  X \s^{-1} Y \ot \s^{-1}  -
  Y \s^{-1} X \ot \s^{-1} + \s^{-1} X \ot \s^{-1} \, = \\
 & - [ Y, \s^{-1} X] \ot \s^{-1} + \s^{-1} X \ot \s^{-1}
 \, = \, 0 \, ,
\end{eqnarray*}
one obtains
 \begin{equation*}
 \tau_2 (TF^\dag) \, = \,
 ( \s^{-1} X - \s^{-1} \d_1 Y ) \ot \s^{-1} Y -
 Y \ot \s^{-1} X \, = \ TF^\dag \, .
\end{equation*}
\end{proof}
\medskip

We conclude with two comments.
The first is that the two periodic Hopf cyclic classes
$TF^\dag$ and $GV^\dag$ are surely nontrivial, because they lift nontrivial
classes. In $\S$\ref{S4.1} it will be shown that they generate
$HP_{(\d , \s^{-1})}^{\ev} (\Hc^\dag_1)$ and, respectively,
$HP_{(\d , \s^{-1})}^{\odd} (\Hc^\dag_1)$.

The second is that not just $\, (\d, \s^{-1} )$ but actually
 every pair $\, (\d, \s^k )$, with $k \in \Zb$,
 is a modular pair in involution for $\Hc_1^{\dagger}$.
Indeed, the twisted antipode
$\, S_\d= \d \ast S$ is an anti-homomorphism that
assumes on the generators the following values:
\begin{eqnarray*}
 S_\d(\s) \, &=&\,   \s^{-1} \, , \qquad \qquad \qquad
 S_\d(Y) \, =\,   -Y + 1 \, , \\
 S_\d(X)  \, &=&\, - \s^{-1} \, (X - \d_1 Y) \, , \quad
  S_\d(\d_1)  \, =\, - \s^{-1} \,  \d_1 \, .
\end{eqnarray*}
One can then easily check that $\, S_\d^2 = \Id$; on
the other hand $\Ad \s^k = \Id$, since $\s$ is central. Again in $\S$\ref{S4.1},
it will be proved that  $HP_{(\d , \s^k)}^* (\Hc^\dag_1) = 0$
for all  $k \neq -1$, $k \in \Zb$.

\subsection{Actions of $\Hc^\dag_1$ on modular Hecke algebras }  \label{S1.3}

The actions of $\Hc^\dag_1$ described in the preceding subsection
automatically factor through the double cover $\Hc^{\dag|2}_1$.
We shall now produce
examples of effective actions of the infinite cyclic cover $\Hc^\dag_1$
on algebras that arise naturally in the theory of modular forms.
\medskip

Let us recall from \cite{cm5} the definition of modular Hecke algebras.
These are algebras
associated to congruence subgroups $\G$ of $\SL(2, \mathbb{Z})$, that
arise from the fusion of the two
quintessential structures coexisting
on modular forms, namely the algebra structure given by the pointwise product
on the one hand, and the action of the Hecke operators on the other. From the
perspective of noncommutative geometry, these algebras describe the
holomorphic `coordinates' of certain noncommutative `arithmetic spaces',
and the action of ${\mathcal H}_1$  is `elliptic', in the sense that it spans the
 `holomorphic tangent space' of the underlying noncommutative space.
\medskip

Referring to \cite{cm5} for the occasionally unexplained term in what follows,
we denote by ${\mathcal M} $
the algebra of holomorphic modular forms of all levels, on which
$G^+(\Qb) := \GL^+(2, \mathbb{\Qb})$ acts via the `slash' operator
\begin{equation*} \label{slash}
f|_{k}\, g \, (z) \, = \, \det (g)^{\frac{k}{2}} \, (cz + d)^{-k} \,
f (g \cdot z) \, ,  \qquad
g = \begin{pmatrix}a  &b\\ c   &d  \end{pmatrix} \in G^+(\Qb) \, .
\end{equation*}
The following variant of the
`slash' operator is more natural when modular forms are
viewed as lattice functions:
 \begin{equation*} \label{dagger}
f\dagger_{k}\, g \, (z) \, =  \, (cz + d)^{-k} \,
f (g \cdot z) \,   = \, (\det g)^{-\frac{k}{2}} \,
f\vert_{k}\, g  \, .
\end{equation*}
\medskip

Given a congruence subgroup $\G \sbs \G(1) = \SL (2, \Zb)$, the
\textit{modular Hecke algebra of level} $\G$, to be denoted
${\mathcal A}(\Gamma)$, is the space of all finitely supported maps
$$
F: \Gamma \backslash G^+(\mathbb{Q}) \rightarrow {\mathcal M} \, ,
\qquad \G \a \mapsto F_{\a}
$$
such that
\begin{equation} \label{tcov}
 F_{\alpha \gamma} \, =  \,
 F_{\alpha} \vert \gamma \, , \quad  \forall \, \alpha \in G^+(\mathbb{Q})
 \, , \quad \forall \, \gamma \in \Gamma \, ,
\end{equation}
endowed with the product
 \begin{equation} \label{conv}
(F^1 *  F^2)_{\alpha} \,:= \,
\sum_{\Gamma \beta \in \Gamma \backslash G^+ (\mathbb{Q})}
\, F^2_{\alpha \beta^{-1}} \vert \beta \cdot F^1_{\beta} \, .
\end{equation}

Any finitely supported map from $\G\bash G^+ (\mathbb{Q})/\G$ to
the algebra $\Mc (\G)$ of modular forms of level $\G$
trivially fulfills the equation (\ref{tcov}), but such maps do not exhaust
all its solutions. In fact, given $f \in \Mc$ there exists an $F \in \Ac(\G)$
such that $\, F_{\a} \, = \, f \, $ if and only if
$\, f \vert \g = f \, , \quad \fl \, \g \in \G \, \cap \a^{-1} \G \,\a$.
 \medskip

To describe the way $\, \Hc_1 \,$ acts
on modular Hecke algebras (cf. \cite[\S 2]{cm2}), we recall that while
the space of  modular forms ${\mathcal M}$ is not invariant under
differentiation, there is a classical operator that corrects the
derivative, namely
\begin{equation*}
X \, = \, \frac{1}{2  \pi i} \frac{d}{dz} - \frac{1}{12  \pi i} \,
\frac{d}{dz} (\log \D) \, Y = \frac{1}{2  \pi i} \,
\frac{d}{dz} - \frac{1}{2  \pi i} \, \frac{d}{dz} (\log \eta^4)
\, Y \, ,
\end{equation*}
where $\, Y \,$ is the grading (by the half-weight) operator
\begin{equation*} \nonumber
Y(f) \, = \, \frac{k}{2} \cdot f \, , \quad \fl \, f \in \Mc_{k} \, ,
\end{equation*}
$\D$ is the discriminant modular form
and $\eta$ is the Dedekind eta function.
\smallskip

With this choice of a `holomorphic connection',
there is a unique action of  $\Hc_1$ on $\Ac (\G)$ sending
the generators of  $\Hc_1$ to the following linear operators on $\Ac (\G)$:
\begin{equation} \nonumber
 Y(F)_{\a} \, = \, Y(F_{\a}) \, ,  \quad
  X(F)_{\a} \, = \, X(F_{\a}) \, ,  \quad
 \d_n (F)_{\a}  \, = \, \mu_{n, \a} \, F_{\a} \, ,
 \end{equation}
 where
\begin{equation} \nonumber
\mu_{\a} \, (z)  =  \frac{1}{2  \pi i} \, \frac{d}{dz} \log
\frac{\eta^4 | \a}{\eta^4}  \, , \qquad
\mu_{n, \, \a} = X^{n-1}(\mu_{\a}) ,
 \quad \fl \, n \geq 1 \, .
\end{equation}
Endowed with this action, $\Ac (\G)$  is an
$\Hc_1$-module algebra, that is
\begin{equation} \label{hact}
 h(a b ) \, = \, \sum h_{(1)} (a) \, h_{(2)} (b) \, ,
 \quad \fl \, a, b \in \Ac(\G) \, , \qquad \fl \, h \in \Hc_1,
\end{equation}
essentially follows from the identity
(see \cite[Lemma 5]{cm2}):
\begin{equation}  \label{xcom}
 X(F \vert {\a})\, = \, X(F)  \vert \a \, + \, \mu_{\a} \cdot Y(F) \vert \a \, , \qquad
  \a \in G^+ (\mathbb{Q}) \, , \quad f \in \Mc_k .
 \end{equation}
Indeed, using  \eqref{xcom} in the third line below, one obtains
\begin{eqnarray} \label{xex}
&\,& X(F^1 \ast F^2)_\a \, = \,   \sum_{\b \in \G\backslash G^+(\Qb)} \, X (
F^2_{\a \b^{-1}} \vert \b  \cdot F^1_\b) \\ \nonumber
&=&\sum_{\b} \, \big(X ( F^2_{\a \b^{-1}} \vert \b ) \cdot F^1_\b
+ F^2_{\a \b^{-1}} \vert \b  \cdot X (F^1_\b) \big) \\ \nonumber
&=&\sum_{\b} \Big( X ( F^2_{\a \b^{-1}} )\vert \b  \cdot F^1_\b
+ \mu_{\b} \, Y  ( F^2_{\a \b^{-1}} ) \vert \b  \cdot F^1_\b +
F^2_{\a \b^{-1}} \vert \b  \cdot X (F^1_\b) \Big) ,
\end{eqnarray}
that is
 \begin{equation*}  
X(F^1 \ast F^2) \, = \, X(F^1)\ast F^2 + F^1 \ast X( F^2) + \d_1(F^1) \ast Y (F^2) \, .
\end{equation*}
In the case of $\d_1$ one uses the $1$-cocycle property of $\mu_{\a}$
(that follows from its very definition) to express
\begin{eqnarray} \label{d1ex}
\d_1(F^1 \ast F^2)_\a &=&  \mu_{\a}
\sum_{\b \in \G\backslash G^+(\Qb) } \,
F^2_{\a \b^{-1}} \vert \b  \cdot F^1_\b \\ \nonumber
&=& ( \mu_{\a \b^{-1}} \vert \b +  \mu_{\b} )
 \sum_{\b} \,  F^2_{\a \b^{-1}} \vert \b ) \cdot F^1_{\b}    \\ \nonumber
&=&  \sum_{\b} \, \big(   ( \mu_{\a \b^{-1}}  F^2_{\a \b^{-1}} )\vert \b  \cdot F^1_{\b}
+ F^2_{\a \b^{-1}}  \vert \b  \cdot  \mu_{\b} \, F^1_{\b} \big),
\end{eqnarray}
which means that $\d_1$ acts as a derivation:
 \begin{equation} \label{d1ex2}
\d_1(F^1 \ast F^2)\, =\, \d_1(F^1) \ast F^2 \, + \, F^1 \ast \d_1( F^2)    \, .
\end{equation}
\medskip

We now define a variant of the algebra $\Ac (\G)$, to be
 denoted $\Ac^{\dagger}(\G)$, obtained by
 modifying the multiplication rule \eqref{conv}  as follows:
 \begin{equation} \label{dconv}
(F^1 *  F^2)_{\alpha} \,:= \,
\sum_{\Gamma \beta \in \G \backslash G^+ (\mathbb{Q})}
\, F^2_{\a \b^{-1}} \dagger \beta \cdot F^1_{\beta} \, ,
\end{equation}

The algebra $\Ac^{\dagger}(\G)$ carries
a natural Hopf module-algebra structure, with respect to
the following `twisted' version of the Hopf algebra $\Hc_1$.
\medskip

 \begin{proposition} \label{thopf}
There exists a unique action of  $\Hc^\dag_1$ on $\Ac^\dag (\G)$ by which
the generators of  $\Hc_1$ act as follows:
\begin{eqnarray} \label{tgen1}
\s (F)_{\a} &=& \det \a \cdot F_{\a} \, , \qquad \fl \, \,
\G \a \in \G\backslash G^+(\Qb) \, ,\\ \label{tgen2}
 Y(F)_{\a}&=&Y(F_{\a}) ,  \quad
  X(F)_{\a}= X(F_{\a}) ,  \quad
 \d_n (F)_{\a} = \mu_{n, \a} \, F_{\a} .
 \end{eqnarray}
 Endowed with this action, $\Ac^\dag (\G)$  is a left
 $\Hc^\dag_1$-module algebra.
 \end{proposition}
\smallskip

\begin{proof} The required verifications are similar to those performed
in the case of $\Hc_1$ (cf.~\cite{cm2}). To illustrate the verification
of the Hopf action property for $X$, we remark that
when the `slash' is replaced with the `dagger' operation, taking into account
that $X(f) $ has weight $k+2$, the identity \eqref{xcom} becomes
\begin{equation} \label{xdcom}
 X(F \dagger {\a})  \, =\,
     \det \a \cdot X(F)  \dagger \a \, + \, \mu_{\a} \cdot Y(F) \dagger \a  \, .
 \end{equation}
Using
 \eqref{xdcom} instead of  \eqref{xcom} in the calculation  \eqref{xex},
one obtains
\begin{eqnarray*}
X(F^1 \ast F^2)_{\a}  &=&  \sum_{\b \in \G\backslash G^+(\Qb)} \, X (
F^2_{\a \b^{-1}} \dagger \b  \cdot F^1_{\b} ) \\
&=&  \sum_{\b} \, \big(X ( F^2_{\a \b^{-1}} \dagger \b ) \cdot F^1_{\b}
+ F^2_{\a \b^{-1}} \dagger \b  \cdot X (F^1_{\b} ) \big) \\
&=&  \sum_{\b} \, \big(
\det \b \, X ( F^2_{\a \b^{-1}} )\dagger \b  \cdot F^1_{\b}  \\
&+&\mu_{\b} \, Y  ( F^2_{\a \b^{-1}} ) \dagger \b  \cdot F^1_{\b}
+ F^2_{\a \b^{-1}} \dagger \b  \cdot X (F^1_{\b} ) \big)  ,
\end{eqnarray*}
whence
\begin{equation}  \nonumber
X(F^1 \ast F^2) \, =\,
 X(F^1)\ast F^2  + \s (F^1)  \ast X( F^2) +
 \d_1(F^1) \ast Y (F^2) \, .
\end{equation}

Also, for $\d_1$, the calculation \eqref{d1ex} takes the form
\begin{eqnarray*}
\d_1(F^1 \ast F^2)_{\a}&=&\mu_{\a}
\sum_{\b \in \G\backslash G^+(\Qb)} \,
F^2_{\a \b^{-1}} \dagger \b  \cdot F^1_{\b}  \\ \nonumber
&=&\big( \det \b\, \mu_{\a \b^{-1}} \dagger \b +  \mu_{\b} \big)
 \sum_{\b} \,  F^2_{\a \b^{-1}} \dagger \b  \cdot F^1_{\b}    \\ \nonumber
&=&\sum_{\b} \, \big(   ( \mu_{\a \b^{-1}}  F^2_{\a \b^{-1}} )\dagger \b
\cdot \det \b \, F^1_{\b}  +
F^2_{\a \b^{-1}}  \dagger \b  \cdot  \mu_{\b} \, F^1_{\b}\big) \, ,
\end{eqnarray*}
or equivalently,
\begin{equation} \nonumber
\d_1(F^1 \ast F^2) =
 \d_1(F^1)\ast F^2  + \s (F^1)  \ast \d_1( F^2)  .
\end{equation}
\end{proof}

\section{Bicocyclic module associated to a
 bicrossed product}\label{S2}
\subsection{MPI coefficients } \label{S2.1}

The Hopf algebra $\Hc_1^\dag$ introduced in $\S$\ref{S1.2}
is isomorphic to the `straight' tensor product
$ \Hc_1\ot \bbc[\mathbb{Z}] $ as an algebra but, as we shall see
in the next section, it  is isomorphic to a cocrossed product
 $\Hc\rt\bbc[\mathbb{Z}]$ as a coalgebra.
In this section we shall associate a bicocyclic module
to such a bicrossed product
Hopf  algebra, whose diagonal can be identified
to the standard cocyclic module defining the Hopf cyclic cohomology.
In view of the 
Getzler-Jones analogue~\cite{gj} of the Eilenberg-Zilber theorem for
bi-paracyclic modules, this makes it possible to compute
the desired Hopf cyclic cohomology out of
the total complex. This method was used by Getzler and Jones to recover
a spectral sequence of Fe\u\i gin and Tsygan~\cite{ft} for 
the cyclic homology of
crossed products of algebras by discrete groups, and was later
generalized by Akbarpour and Khalkhali~ \cite{ak1} to the the cyclic homology of
crossed products of algebras by Hopf algebras.
The key obstacle in implementing
this type of construction for
Hopf cyclic cohomology resides in the
interaction between the antipode and coaction, which 
is tackled in Lemma \ref{antipode} below. The pay-off for the extra 
difficulty is that we obtain a bicocyclic module, not just a cylindrical one.
 
\bigskip
We start with
a Hopf algebra $\Hc$ which admits a left coaction by a commutative Hopf algebra $\Kc$,
$ \, \rho: \Hc\ra \Kc\ot \Hc$,
for which we use the (Sweedler-type) notation
\begin{equation*}
\rho(h)=h\ns{-1}\ot h\ns{0},\quad \rho^{(k)}(h)=h\ns{-k}\odots h\ns{-1}\ot h\ns{0},
\end{equation*}
where $\quad \rho^{(k)}=\Id_{\Kc}\ot \rho^{(k-1)} \quad$ and $\quad\rho^{(1)}=\rho$.

 We call $\Hc$ a $\Kc$-\textit{comodule algebra}
 if for any $h,g\in \Hc$ the following hold:
 \begin{align}
 &\rho(hg)= h\ns{-1}g\ns{-1}\ot h\ns{0}g\ns{0}\label{c-a-1} \, ,\\
 &\rho(1)=1\ot 1\label{c-a-2}.
 \end{align}
Similarly, we call  $\Hc$ a $\Kc$-\textit{comodule coalgebra}  if for any $h\in\Hc$
 \begin{align}
 &h\ps{1}\ns{-1}h\ps{2}\ns{-1}\ot h\ps{1}\ns{0}\ot h
 \ps{2}\ns{0}= h\ns{-1}\ot h\ns{0}\ps{1}\ot h\ns{0}\ps{2},\label{c-c-1}\\
 &\epsilon(h\ns{-1})\ot h\ns{0}=1\ot h\label{c-c-2}.
 \end{align}

 \begin{definition} Let  $\Hc$  be a Hopf algebra which
admits  a left coaction by a Hopf algebra $\Kc$. If
 via this coaction $\Hc$ is simultaneously a $\Kc$-comodule algebra and
 a $\Kc$-comodule coalgebra, then
$\Hc$ will be called a $\Kc$-{\em comodule Hopf algebra}.
\end{definition}

\bigskip We recall that if   $\Hc$ is a left $\Kc$-comodule coalgebra
then the cocrossed product coalgebra  $\Hc\rt \Kc$  has  $\Hc\ot \Kc$
as underlying vector space and the following coalgebra structure:
\begin{align}
&\Delta(h\rt k)= h\ps{1}\rt h\ps{2}\ns{-1}k\ps{1}\ot h\ps{2}\ns{0}\rt k\ps{2}, \\
&\epsilon(h\rt k)=\epsilon(h)\epsilon(k).
\end{align}
\begin{lemma}\label{crossed}
The cocrossed product coalgebra  $\Hc\rt \Kc$ admits a
Hopf algebra structure via the tensor
product algebra structure and with the antipode given by
\begin{equation*}
S(h\rtimes k)=S(h\ns{0})\rtimes S(h\ns{-1}k).
\end{equation*}
\end{lemma}

\begin{proof}
We need to check that $\Delta$ and $\epsilon$ are
  algebra map and  that $S$ is convolution inverse of $\Id$.
  In order to show that $\Delta$ is an algebra map we will use the
  commutativity of $\Kc$ and also the fact that $\Hc$ is
 $\Kc$-comodule algebra, i.e. (\ref{c-a-1}) and (\ref{c-a-2}). Thus,
\begin{align*}
&\Delta((h^1\rt k^1)(h^2\rt k^2))=\Delta(h^1h^2\rt k^1k^2)=\\
&=(h^1\ps{1}h^2\ps{1})\rt
(h^1\ps{2}h^2\ps{2})\ns{-1}k^1\ps{1}k^2\ps{1} \ot
(h^1\ps{2}h^2\ps{2})\ns{0}\rt k^1\ps{2}k^2\ps{2}= \\
&=(h^1\ps{1}\rt h^2\ps{2}\ns{-1}k^1\ps{1}\ot h^1\ps{2}\ns{0}\rt
k^1\ps{2})\\
&(h^2\ps{1}\rt h^1\ps{2}\ns{-1}k^2\ps{1}\ot h^2\ps{2}\ns{0}\rt
k^2\ps{2})=\Delta(h^1\rt k^1)\Delta(h^2\rt k^2).
\end{align*}
Similarly, by relying on the fact that $\epsilon_{\Hc}$ and
$\epsilon_{\Kc}$ are algebra maps, one easily checks that  $\epsilon$ is algebra
maps too.

To verify the antipode axioms, we write
\begin{align*}
&S\ast \Id(h\rt k)=S(h\ps{1}\rt
h\ps{2}\ns{-1}k\ps{1})(h\ps{2}\ns{0}\rt k\ps{2})\\
&=S(h\ps{1}\ns{0})\rt
S(h\ps{1}\ns{-1}h\ps{2}\ns{-1}k\ps{1})(h\ps{2}\ns{0}\rt k\ps{2})\\
&=S(h\ns{0}\ps{1})\rt S(h\ns{-1})(h\ns{0}\ps{2}\rt
k\ps{2})\\
&=S(h\ns{0}\ps{1})h\ns{0}\ps{2}\rt
S(h\ns{-1})S(k\ps{1})k\ps{2}=\epsilon(h\ns{0})1_\Hc\rt
S(h\ns{-1})\epsilon(k)\\
&=\epsilon(h)1_\Hc\rt \epsilon(k)1_\Kc.
\end{align*}
A similar computation shows that
$$\Id\ast S(h\rt
k)=\epsilon(h\rt k)1_\Hc\rt 1_\Kc.
$$
\end{proof}

\begin{lemma}\label{antipode}
If $\Hc$ is a $\Kc$-comodule  Hopf algebra then the antipode of $\Hc$ is $\Kc$-colinear:
\begin{equation*}
\quad h\ns{-1}\ot S(h\ns{0})= S(h)\ns{-1}\ot S(h)\ns{0},\, \, h \in \Hc.
\end{equation*}
\end{lemma}
\begin{proof}
Consider $R_1$, and $R_2$ in $\Hom(\Hc, \Kc\ot \Hc)$, defined
by:
\begin{align*}
&R_1(h):=h\ps{1}\ns{-1}h\ps{2}\ns{-1}\ot  h\ps{1}\ns{0}
S(h\ps{2}\ns{0}) \, ,\\
&R_2(h)=1\ot \epsilon(h).
\end{align*}
We first prove that  $R_1\,  = \, R_2$.
\smallskip

On applying $(1\ot m_\Hc)\circ(1\ot 1\ot S)$, where $m_\Hc$ is the
multiplication of $\Hc$, to both sides of
 (\ref{c-c-1}) one obtains:
 \begin{align*}
&R_1(h)=h\ps{1}\ns{-1}h\ps{2}\ns{-1}\ot  h\ps{1}\ns{0}
S(h\ps{2}\ns{0})=\\
& h\ns{-1}\ot h\ns{0}\ps{1}S(h\ns{0}\ps{2})=\\
&h\ns{-1}\ot \epsilon(h\ns{0})= 1\ot\epsilon(h)=R_2(h)
\end{align*}

To complete the proof it suffices to convolute $R_1$ and $R_2$ with the operator
$L \in \Hom(\Hc, \Kc\ot \Hc)$,
\begin{equation*}
L(h):=S(h)\ns{-1}\ot S(h)\ns{0}.
\end{equation*}

 Indeed, one has
 \begin{equation*}
   h\ns{-1}\ot S(h\ns{0})= L\ast R_1(h)=L\ast
R_2(h)=S(h)\ns{-1}\ot S(h)\ns{0}.
\end{equation*}
\end{proof}

\begin{lemma}\label{twisted}
If $\alpha:\Hc\rightarrow \bbc$ is a $\Kc$ colinear character, then
so is $S_\alpha=\alpha\ast S$.
\end{lemma}
\smallskip

\begin{proof}  On the one hand,
\begin{eqnarray*}
 S_\alpha(h)\ns{-1}\ot
S_\alpha(h)\ns{-0}&=&\alpha(h\ps{1})S(h\ps{2})\ns{-1}\ot
S(h\ps{2})\ns{0}\\
&=& h\ps{2}\ns{-1}\ot \alpha(h\ps{1})S(h\ps{2}\ns{0}) \, ,
\end{eqnarray*}
since, by Lemma \ref{antipode}, $S$ is colinear.

 On the other hand,
\begin{eqnarray*}
 h\ns{-1}\ot S_\alpha(h\ns{0})&=&h\ns{-1}\ot \alpha(h\ns{0}\ps{1})
 S(h\ns{0}\ps{2})\\
 &=&h\ps{1}\ns{-1}h\ps{2}\ns{-1}\ot
 \alpha(h\ps{1}\ns{0})S(h\ps{2}\ns{0})\\
 &=&h\ps{2}\ns{-1}\ot \alpha(h\ps{1})S(h\ps{2}\ns{0}),
 \end{eqnarray*}
 using the colinearity of $\alpha$.
\end{proof}
\medskip

\begin{lemma}\label{tantipode}
Let $\alpha$ and $\beta$ be characters for $\Hc$ and $\Kc$
respectively. If $\alpha$ is $\Kc$ colinear, then
$$S_{\alpha\ot \beta}(h\rtimes k)=S_\alpha(h\ns{0})\rtimes S_\beta(h\ns{-1}k).$$
\end{lemma}

\begin{proof}
\begin{align*}
&S_{\alpha\ot \beta}(h\rtimes k)=(\alpha\ot \beta)(h\ps{1}\rtimes
h\ps{2}\ns{-1}k\ps{1})S(h\ps{2}\ns{0}\rtimes k\ps{2})\\
&=\alpha(h\ps{1})S(h\ps{2}\ns{0})\rtimes
\beta(h\ps{2}\ns{-2}k\ps{1}) S(h\ps{2}\ns{-1}k\ps{2})\\
&=\alpha(h\ps{1})S(h\ps{2}\ns{0})\rtimes
S_\beta(h\ps{2}\ns{-1})S_\beta(k).
\end{align*}
On the other hand,
\begin{align*}
&S_\alpha(h\ns{0})\rtimes
S_\beta(h\ns{-1}k)=\alpha(h\ns{0}\ps{1})S(h\ns{0}\ps{2})\rtimes
S_\beta(h\ns{-1})S_\beta(k)\\
&=\alpha(h\ps{1}\ns{0})S(h\ps{2}\ns{0})\rtimes
S_\beta(h\ps{1}\ns{-1}h\ps{2}\ns{-1})S_\beta(k)\\
&=\alpha(h\ps{1})S(h\ps{2}\ns{0})\rtimes
S_\beta(h\ps{2}\ns{-1})S_\beta(k),
\end{align*}
where in the last equality we used once more the colinearity of $\alpha$.
\end{proof}
\medskip

\begin{definition}
Let $\Hc$ be a  $\Kc$-comodule Hopf algebra. An MPI
$(\alpha,\mu)$ over $\Hc$ will be called $\Kc$-{\em coinvariant} if $\alpha$ is
colinear and $\mu\ns{-1}\ot \mu\ns{0}=1\ot \mu$.

A character $\beta$ of $\Kc$ will be called
 {\em stable} if $\beta(h\ns{-1})h\ns{0}=h$,  for all $h\in \Hc$.
\end{definition}
\medskip

\begin{proposition}\label{MPI}
Let $(\alpha, \mu)$  and $(\beta, \nu)$ be MPIs for $\Hc$ and $\Kc$
respectively. If $(\alpha,\mu)$ is coinvariant and $\beta$ is
stable, then the pair $(\alpha\ot \beta,\mu\rtimes\nu)$ is an MPI for
$\Hc\rtimes \Kc$.
\end{proposition}

\begin{proof}  We express
\begin{eqnarray*}
S_{\alpha\ot\beta}^2(h\rtimes k)&=& S_{\alpha\ot\beta}(S_\alpha(h\ns{0})\rtimes
S_\beta(h\ns{-1}k))\\
&=&S_\alpha(S_\alpha(h\ns{0})\ns{0})\rtimes
S_\beta(S_\alpha(h\ns{0})\ns{-1}S_\beta(h\ns{-1}k)), \\
&=&S_\alpha^2(h\ns{0})\rtimes
S_\beta(h\ns{-1}S_\beta(h\ns{-2}))S_\beta^2(k)\\
&=&\mu h\ns{0}\mu^{-1}\rtimes \beta(h\ns{-1})\nu k\nu^{-1}\\
&=&\mu h\mu^{-1}\rt \nu k\nu^{-1},
\end{eqnarray*}
using first the property that  $S_\alpha$ is $\Kc$-colinear, then
the stability of  $\beta$.
\end{proof}

\smallskip
From now on we fix MPIs  $(\alpha, \mu)$ and $(\beta, \nu)$
satisfying the assumptions of the previous proposition, and thus
such that $(\alpha\ot\beta, \mu\rt\nu)$ is an MPI for $\Hc\rt\Kc$.
The rest of this section is devoted to the Hopf cyclic cohomology of
$\Hc\rt \Kc$  with coefficients in  the MPI $(\alpha\ot\beta,
\mu\rt\nu).$

\bigskip

For the reader's convenience,  we recall that
the Hopf cyclic cohomology of a Hopf algebra $\Hc$
 with respect to an MPI  $(\delta,\sigma)$ is  defined as the cyclic cohomology of
  the cocyclic module $\Hc_\nr = \{\Hc^n_\nr \, ; \, n \geq 0 \}$, where
  \begin{equation*}
  \Hc^n_\nr := \Hc^{\ot n} , \quad  \text{if} \, \, n\ge 1 \quad \text{and} \quad
   \Hc^0_\nr := \bbc \, ,
  \end{equation*}
 with the following cocyclic structure:
\begin{align}
&\p_0(h^1\ot\dots\ot h^q)=  1\ot h^1\ot\dots\ot h^q \label{CM1}
\\
&\p_i(h^1\ot\dots\ot h^q)=h^1\ot\dots\ot \Delta(h^i)\ot\dots \ot h^q\\
&\p_{n+1}(h^1\ot\dots\ot h^q)=h^1\ot \dots \ot  h^q\ot \sigma \\
&\sigma_i(h^1\ot\dots\ot h^q)= h^1\ot \dots \ot \epsilon(h^{j+1})\ot\dots\ot h^q\\
&\tau(h^1\ot\dots\ot h^q)= S_\d (h^1)\cdot(h^2\ot\dots\ot h^q\ot
\sigma)\label{CM6}.
\end{align}

In what follows, we shall often employ the following shorthand
notation:
\begin{align*}
&\tilde{h}=h^1\odots  h^q\, , \qquad \qquad
 \tilde{k}=k^1\odots  k^p \, ,\\
 &\td{h}\ns{-1}\ot \td{h}\ns{0}= (h^1\ns{-1}\dots h^q\ns{-1}) \ot (h^1\ns{0}\odots h^q\ns{0}).
 \end{align*}

We now recall that a {\em bicocyclic module} is a bigraded module
whose rows and  columns  are  cocyclic modules and, in addition,
vertical operators  commute with  horizontal operators. It is known,
cf. \cite{gj}, that the diagonal of a  bicocyclic module is a
cocyclic module.

We define  the bicocyclic module $\FC = \{\FC^{(p,q)} \, ; \, (p,
q) \in \Zb^+ \times \Zb^+ \}$ as follows:
\begin{equation}\label{bisim}
 \FC^{(p,q)} := \Kc^{p} \ot \Hc^{q}  \, ,
\end{equation}
the horizontal maps
\begin{align*}
&\hd_i:\FC^{(p,q)}\rightarrow \FC^{(p+1,q)}, && 0\le i\le p+1\\
&\hs_j: \FC^{(p,q)}\rightarrow \FC^{(p-1,q)}, &&0\le j\le p-1\\
&\hta:  \FC^{(p,q)}\rightarrow \FC^{(p,q)},
\end{align*}
are defined by:
\begin{align}\label{hd0}
&\hd_0(\td{k}\ot \td{h})= 1\ot k^1\ot\dots\ot k^p\ot \td{h}\\ \label{hdi}
&\hd_j(\td{k}\ot \td{h})= k^1\ot\dots\ot \Delta(k^i)\ot\dots \ot k^p\ot \td{h}\\ \label{hdn}
&\hd_{p+1}(\td{k}\ot \td{h})=k^1\ot \dots \ot k^p\ot \td{h}\ns{-1}\nu\ot  \td{h}\ns{0}\\ \label{hs}
&\hs_j(\td{k}\ot \td{h})= k^1\ot \dots \ot \epsilon(k^{j+1})\ot\dots\ot k^p\ot \td{h}\\ \label{ht}
&\hta(\td{k}\ot \td{h})= S_\beta(k^1)\cdot(k^2\ot\dots\ot k^p\ot
\nu\td{h}\ns{-1})\ot \td{h}\ns{0};
\end{align}
the vertical structure is just the cocyclic structure of $\Hc$, with
\begin{align} \label{vd}
&\vd_i= \Id\ot \p_i:\FC^{(p,q)}\rightarrow \FC^{(p,q+1)}, &&0\le i\le q+1\\ \label{vs}
&\vs_j=\Id\ot \sigma_j: \FC^{(p,q)}\rightarrow \FC^{(p,q-1)}, && 0\le j\le q-1\\ \label{vt}
&\vta= \Id\ot \tau:  \FC^{(p,q)}\rightarrow \FC^{(p,q)},
\end{align}
where $\p_i$, $\sigma_i$, and $\tau$ are given by (\ref{CM1})-(\ref{CM6}).
\medskip

Recall from  \cite{hkrs1} that {\em a right-left stable anti
Yetter-Drinfeld} (SAYD for short) module over  a Hopf algebra  $\Hc$
is a right $\Hc$-module $M$ equipped with a left coaction of $\Hc$
such that for any $m\in M$ and any $h\in \Hc$,
\begin{eqnarray} \nonumber
m\sns{0}m\sns{-1}&=&m \, , \\ \label{SAYD}
(mh)\sns{-1} \ot (mh)\sns{0} &=& S(h\ps{3})m\sns{-1}h\ps{1}\ot m\sns{0}h\ns{2};
\end{eqnarray}
we have added the bar superscript to the notation as an extra
precaution, to avoid any possible
confusion in case $\Hc$ itself is a comodule for some coaction.
 \medskip

 The Hopf cyclic cohomology of a Hopf algebra $\Hc$ with coefficients in a SAYD module
 $M$ is  defined as the cyclic cohomology of the cocyclic module $C^*_H(\Hc,M): =
  \{\Hc ^{\ot n}\ot M\}_{n\ge 0} \, $ endowed with the following operators:
 \begin{align*}
 &\p_0(h^1\odots h^n\ot m)=1\ot h^1\ot h^1\odots h^n\ot m\\
 &\p_i(h^1\odots h^n\ot m)=h^1\odots\Delta(h^i)\odots h^n\ot m\\
 &\p_{n+1}(h^1\odots h^n\ot m)= h^1\odots h^n \ot m\sns{-1}\ot m\sns{0} \\
 &\sigma_i(h^1\odots h^n\ot m)=h^1\odots \epsilon(h^{i+1})\odots h^n\ot m\\
 &\tau(h^1\odots h^n\ot m)=S(h^1\ps{2})\cdot h^2\odots h^n\ot m\sns{-1}\ot m\sns{0}h\ps{1}.
 \end{align*}
\medskip

We now equip each $\Hc^{\ot q}$ with a right module and left
comodule structure over  $\Kc$ via $\, \, \td hk=\beta(k)\td h $,
respectively $\bar \rho: \Hc^{\ot q}\ra \Kc\ot \Hc^{\ot q}$ defined by
$$\bar\rho(\td h)=\td h\ns{-\bar 1}\ot \td h\ns{\bar 0}:= \td h\ns{-1}\nu\ot\td h\ns{0} .
$$
\smallskip

\begin{lemma} \label{lSAYD} For each $q\ge 0$, $\Hc^{\ot q}$ is an  SAYD module over
$\Kc$ and the $q$th row of \eqref{bisim} is isomorphic to the cocyclic module
$\{C^p_{\Kc}(\Kc,\Hc^{q})\}_{p\ge 0}.$
\end{lemma}
\smallskip

\begin{proof}
The second part of the lemma is clear, so
 we only need to prove the first statement.
 Using the commutativity of $\Kc$ and also the fact that $(\delta,\nu)$ is a MPI for $\Kc$,
 one gets
$$\bar\rho(\td hk)=\beta(k)\bar\rho(\td h)=\beta(k)\td{h}\ns{-1}\nu\ot\td{h}\ns{0}=
S(k\ps{3})\td h\ns{-1}\nu k\ps{1}\ot \beta(k\ps{2}) \td h\ns{0}, $$
which shows that $\Hc^{\ot q}$ is a anti Yetter-Drinfeld module over
$\Kc$. To prove stability we use the extra property of $\beta$, i.e.
$\beta(h\ns{-1})h\ns{0}=h$.
\end{proof}
\medskip

\begin{proposition}
The bigraded module $\FC$ defined by \eqref{bisim} is bicocyclic.
\end{proposition}
\smallskip

\begin{proof}
First we show that each row and column is a cocyclic module. For
columns this is obvious because the operators are identical to those
of the cocyclic module $\Hc_\nr$. On the other hand, Lemma
\ref{lSAYD} shows that the rows are cocyclic module too.

It remains to show that  the vertical operators and the horizontal
operators commute. To this end, it suffices  to prove that $\hta$
commutes with all vertical operators. We only show that $\hta$ and
$\vta$ commute with one another and leave to the reader to check the
rest. Using Lemma \ref{antipode}, Lemma  \ref{twisted}, the fact
that $\mu$ is $\Kc$-coinvariant, (\ref{c-a-1}) and  (\ref{c-c-1}),
we can write:
\begin{align*}
&\hta\circ\vta(\td k\ot \td h)= \hta(\td k\ot S_\alpha(h^1)\cdot h^2\odots h^q\ot \mu)\\
&=S_\beta(k^1)\cdot(k^2\odots k^p\ot \nu
S(h^1\ps{q+1})\ns{-1}h^2\ns{-1}
\dots  \\
&\dots S(h^1\ps{2})\ns{-1} h^q\ns{-1}\ot S_\alpha(h^1\ps{1})\ns{-1}\mu\ns{-1})\ot S(h^1\ps{q+1})\ns{0}h^2\ns{0}\ot \dots\\
&\dots\ot  S(h^1\ps{2})\ns{0} h^q\ns{0}S_\alpha(h^1\ps{1})\ns{0}\mu\ns{0}=\\
&=S_\beta(k^1)\cdot(k^2\odots k^p\ot \nu \td h\ns{-1})\ot S_\alpha(h^1\ns{0})\cdot(h^2\ns{0}\odots h^q\ns{0}\ot \mu)\\
&=\vta\circ\hta(\td k\ot \td h).
\end{align*}
\end{proof}
\medskip

 Let $\FD$ be the
diagonal of $\FC$, that is $\FD= \{\FD^n \, ; \, n \geq 0 \}$ with
  \begin{equation*}
 \FD^n:= \Kc^{\ot n}\ot \Hc^{\ot n} \equiv \Kc^{\ot n}\ot   \Hc^{\ot n},
  \quad  \text{if} \, \, n\ge 1 \quad \text{and} \quad
   \FD^0 := \bbc \, ,
  \end{equation*}
equipped with the following cocyclic structure:
\begin{align*}
&d_0(\tilde{k}\ot   \tilde{h})= 1\ot k^1\ot\dots \ot k^n\ot   1\ot h^1\ot \dots \ot h^n, \\
&d_i(\tilde{k}\ot   \tilde{h})=k^1\ot\dots\ot
\Delta(k^i)\ot\dots \ot k^n\ot   h^1\dots \ot  \Delta(h^i) \ot \dots\ot h^n,\\
&d_{n+1}(\tilde{k}\ot   \tilde{h})= k^1\ot\dots
\ot k^n\ot \nu\tilde{h}\ns{-1}\ot    h^1\ns{0}\ot \dots \ot h^n\ns{0}\ot \mu,\\
&s_j(\tilde{k}\ot   \tilde{h})= k^1\ot\dots\ot
\epsilon(k^{j+1})\ot\dots \ot k^n\ot   h^1
\dots \ot    \epsilon(h^{j+1}) \ot \dots\ot h^n,\\
&t(\tilde{k}\ot   \tilde{h})= S_\beta(k^1)\cdot (k^2\ot\dots\\
&\dots \ot k^n\ot \nu \tilde{h}\ns{-1}) \ot  
S_\alpha(h^1\ns{0})\cdot(h^2\ns{0}\ot \dots\ot h^n\ns{0}\ot \mu).
\end{align*}

We now proceed to identify the cocyclic module
$(\Hc\rt \Kc)_\nr$ with $\FD$.
To this end we introduce the maps $\Psi$ and its inverse  $\Psi^{-1}$.
First, we define

$\, \Psi: (\Hc\rt \Kc)_\nr\longrightarrow \FD$, by setting
\begin{equation}\label{PSi}
\Psi= \prod_{i=1}^{n}  \Id_{\Kc^{\ot n-i}}\ot \top^{\ot i}\ot  \Id_{\Hc^{\ot n-i}},
\end{equation}
where
\begin{align}
\top:\Kc\ot \Hc \ra \Kc\ot \Hc, \quad \top(h\ot k)= h\ns{-1}k\ot h\ns{0}.
\end{align}
The effect of $\Psi$ is explicitly seen as follows:
 \begin{align*}
&\Psi(h^1\rt k^1\ot \dots\ot h^n\rt k^n)=
 h^1\ns{-n}k^1 \ot h^1\ns{-n+1}h^2\ns{-n+1}k^2\ot\dots\\
&\dots\ot h^1\ns{-2}\dots h^{n-1}\ns{-2}k^{n-1}
\ot  h^1\ns{-1}\dots h^n\ns{-1}k^n\ot   h^1\ns{0}\odots
h^n\ns{0}.
\end{align*}
The inverse  $\, \Psi^{-1}:  \FD\longrightarrow  (\Hc\rt \Kc)_\nr $ is
given by a similar expression,
 \begin{equation}\label{PSi-1}
\Psi^{-1}= \prod_{i=0}^{n-1}  id_{\Kc^{\ot i}}\ot \bot^{\ot n-i}\ot  id_{\Hc^{\ot i}},
\end{equation}
where
\begin{align}
\bot:\Kc\ot \Hc \ra  \Hc\ot \Kc, \quad \bot(k\ot h)= h\ns{0}\ot S(h\ns{-1})k.
\end{align}
Its explicit expression is:
 \begin{align*}
&\Psi^{-1}(\tilde{k}\ot   \tilde{h})=\\
&h^1\ns{0}\rt S(h^1\ns{-1})
k^1\ot\dots \ot h^n\ns{0}\rt S(h^1\ns{-n}h^2\ns{-n+1}\dots
h^n\ns{-1})k^n .
\end{align*}
In order to check that
$\Psi$ and $\Psi^{-1}$ are in fact inverse of one another,
 we only need to show  that $\top$ and $\bot$  are so.
Indeed,
\begin{align*}
&\top\bot(k\ot h)=\top(h\ns{0}\ot S(h\ns{-1})k)= h\ns{-1}S(h\ns{-2})k\ot h\ns{0}=k\ot h. \\
&\bot\top(h\ot k)=\bot(h\ns{-1}k\ot h\ns{0})= h\ns{0}\ot S(h\ns{-1})h\ns{-2}k= h\ot k. \\
\end{align*}
\begin{proposition} \label{map}
The maps $\Psi $ and $\Psi^{-1}$  are cyclic.
\end{proposition}
\begin{proof}
 To show that $\Psi$ is a cyclic map,  we first use Lemma
\ref{tantipode} to observe that:
\begin{align*}
&\Delta^n(S_{(\alpha\ot \beta)}(h\ns{0}\rt S(h\ns{-1}k)))=
S(h\ps{n+1})\rt h\ps{n}\ns{-1}\dots h\ps{1}\ns{-n}S(k\ps{n+1})\ot \\
&\ot S(h\ps{n}\ns{0})\rt h\ps{n-1}\ns{-1}\dots h\ps{1}\ns{-n+1}
S(k\ps{n})\ot\dots \ot  S(h\ps{2}\ns{0})\rt h\ps{1}\ns{-1}S(k\ps{2})\\
&\ot S_\alpha(h\ps{1}\ns{0})\rt S_\beta(k\ps{1}).
\end{align*}
In view of the above identity we can write
\begin{align*}
&\tau\Psi^{-1}(\td{k}\ot  \td h)=\\
& =\tau(h^1\ns{0}\rt S(h^1\ns{-1})
k^1\ot\dots \ot h^n\ns{0}\rt S(h^1\ns{-n}h^2\ns{-n+1}\dots h^n\ns{-1})k^n )\\
&=S_{\alpha\ot\beta}(h^1\ns{0}\rt S(h^1\ns{-1}k^1)\cdot( h^2\ns{0}\rt S(h^1\ns{-2}h^2\ns{-1})k^2\ot\dots\\
&\dots\ot  h^n\ns{0}\rt S(h^1\ns{-n}h^2\ns{-n+1}\dots h^n\ns{-1})k^n\ot  \mu\rt\nu) =\\
&=\Delta(S_{\alpha\ot\beta}(h^1\ns{0}\rt S(h^1\ns{-1}k^1))( h^2\ns{0}\rt S(h^1\ns{-2}h^2\ns{-1})k^2\ot\dots\\
&\dots\ot  h^n\ns{0}\rt S(h^1\ns{-n}h^2\ns{-n+1}\dots h^n\ns{-1})k^n\ot  \mu\rt\nu) =\\
&=S(h^1\ns{0}\ps{n+1})h^2\ns{0} \rt h^1\ns{0}\ps{n}\ns{-1}\dots h^1\ns{0}\ps{1}
\ns{-n}S(k^1\ps{n+1})S(h^1\ns{-2}h^2\ns{-1})\ot \dots\\
&\dots\ot S(h^1\ns{0}\ps{2}\ns{0}) h^n\ns{0}\rt h^1\ns{0}\ps{1}\ns{-1}S(k^1
\ps{2})S(h^1\ns{-n}h^2\ns{-n+1}\dots h^n\ns{-1})k^n\ot\\
&\ot S_\alpha(h^1\ns{0}\ps{1}\ns{0})\mu\rt S_\beta(k^1\ps{1})\nu.
\end{align*}
Using  the fact that $\Delta^{(n)}: \Hc\ra \Hc^{\ot n+1}$,   $S$,
and $S_{\alpha}$ are $\Kc$-colinear,  and $\Kc$ is commutative,  one has:
\begin{align*}
&\Psi\tau\Psi^{-1}(\td k\ot   \td h)=
h^1\ns{0}\ps{n+1}\ns{-n}\dots
h^1\ns{0}\ps{1}\ns{-n}h^2\ns{0}\ns{-n} \\
&S(h^1\ns{-1}h^2\ns{-1})S(k^1\ps{n+1})k^2
\ot h^1\ns{0}\ps{n+1}\ns{-n+1}\dots\\
&\dots h^1\ns{0}\ps{1}\ns{-n+1}h^2\ns{0}\ns{-n+1}h^3\ns{0}\ns{-n+1}\\
&S(h^1\ns{-2}
h^2\ns{-2}h^3\ns{-1})S(k^1\ps{n})k^3\ot\\
&\odots  h^1\ns{0}\ps{n+1}\ns{-2}\dots
h^1\ns{0}\ps{1}\ns{-2}h^2\ns{0}\ns{-2}\dots h^n\ns{0}\ns{-2}\\
&S(h^1\ns{-n+1}
h^2\ns{-n+1}\dots h^n\ns{-1})S(k^1\ps{2})k^n\\
&\ot h^1\ns{0}\ps{n+1}\ns{-1}\dots h^1\ns{0}\ps{1}\ns{-1}h^2\ns{0}\ns{-1}
\dots h^n\ns{0}\ns{-1}S_\beta(k^1\ps{1})\nu\ot   \\
&S(h^1\ns{0}\ps{n+1}\ns{0})h^2\ns{0}\ns{0}\odots S(h^1\ns{0}\ps{1}\ns{0})\mu \\
&=h^1\ns{-n}h^2\ns{-n}S(h^1\ns{-n-1}h^2\ns{-n-1})S(k^1\ps{n+1})k^2\ot
\dots \\
&\dots\ot h^1\ns{-2}\dots
h^n\ns{-2}S(h^1\ns{1-2n}h^2\ns{1-2n}h^3\ns{3-2n}\dots
h^n\ns{-3})S(k^1\ps{2})k^n\\
&\ot h^1\ns{-1}\dots h^n\ns{-1} S_\beta(k^1\ps{1})\nu\ot  
S(h^1\ns{0}\ps{n+1})h^2\ns{0}\ot \dots \\
&\dots\ot S(h^1\ns{0}\ps{2})h^n\ns{0}\ot
S_\alpha(h^1\ns{0}\ps{1})\mu= \\
&= S(k^1\ps{n+1})\odots S(k^1\ps{2})\ot S_\beta(k^1\ps{1})\nu \td
h\ns{-1}\ot   S(h^1\ns{0}\ps{n+1})h^2\ns{0}\ot \\
&\dots\ot S(h^1\ns{0}\ps{2})h^n\ns{0}\ot
S_\alpha(h^1\ns{0}\ps{1})\mu=     t(\td k\ot   \td h).
\end{align*}
\ni By a similar computation one shows that
$\Psi\p_i\Psi^{-1}=d_i$ and $\Psi\sigma_i\Psi^{-1}=s_i$.
\end{proof}
\medskip

Let $Tot(\FC)$  the corresponding total mixed complex
$$Tot(\FC)_n=\bigoplus_{p+q=n} \Kc^{p}\ot \Hc^{q}\, .
$$
We shall denote by $\Tot(\FC)$ the associated \textit{normalized} subcomplex,
obtained by retaining only the elements annihilated by all degeneracy
operators. The total  boundary $b_T+B_T$ defined by
\begin{align}
&\hb_p=\sum_{i=0}^{p+1} (-1)^{i}\hd_i,&&\vb_q=\sum_{i=0}^{q+1} (-1)^{i}\vd_i, \label{b}\\
&b_T=\sum_{p+q=n}\hb_p+\vb_q,\\
&\hB_p=(\sum_{i=0}^{p-1}(-1)^{(p-1)i}\hta^i)\hs_{p-1}\hta,  &&\vB_q=
(\sum_{i=0}^{q-1}(-1)^{(q-1)i}\vta^i)\hs_{q-1}\vta,\\
&B_T=\sum_{p+q=n}\hB_p+\vB_q.\
\end{align}

The total complex of a bicocyclic module is a mixed complex, and so
 its cyclic cohomology is well-defined. In view of the analogue of
 the Eilenberg-Zilber theorem for bi-paracyclic
modules~\cite{gj}, we can assert that the diagonal mixed complex
$(\FD, b, B)$ and
the total mixed complex $(\Tot\FC,b_T, B_T)$ are quasi-isomorphic in
Hochschild and cyclic cohomology, via Alexander-Whitney map

\begin{equation}\label{AW}
\, AW:= \bigoplus_{p+q=n} AW_{p,q}: \Tot(\FC)^n\ra \FD^n,
\end{equation}
\begin{equation*}
AW_{p,q}: \Kc^{\ot p}\ot  \Hc^{\ot q}\longrightarrow \Kc^{\ot
p+q}\ot  \Hc^{\ot p+q}
\end{equation*}
\begin{equation*}
AW_{p,q}=(-1)^{p+q}\underset{
p\,\text{times}}{\underbrace{\vd_0\vd_0\dots
\vd_0}}\hd_n\hd_{n-1}\dots \hd_{p+1} \, .
\end{equation*}

Together with  Proposition \ref{map}, this allows to conclude with the
main result of this section.
\medskip

\begin{theorem}\label{mix}
The mixed complex $((\Hc\rt\Kc)_\nr, b, B)$ is quasi-isomorphic to
the mixed complex $(\Tot\FC,b_T, B_T)$.
\end{theorem}

\subsection{SAYD coefficients}  \label{S2.2}

With a view towards future applications of the methods of the present paper,
we now digress to extend the results of the previous subsection
to the general framework of Hopf cyclic cohomology for
module coalgebras with coefficients in SAYD modules (cf.~\cite{hkrs1,
hkrs2}). The allowable SAYD coefficients are those that
satisfy appropriate compatibility conditions,
matching the cocrossed product construction of the Hopf algebra itself.
\smallskip

With $\Hc$ and $\Kc$ satisfying the same assumptions
as above, we let $C$ and $D$ denote a left $\Hc$-module coalgebra,
respectively a left $ \Kc$-module coalgebra. In addition we assume
that  $C$ is left $\Kc$-comodule coalgebra  via
\begin{align}
&{\cal r} : C\ra \Kc\ot C, \\
&{\cal r}(c)=c\ns{-1}\ot c\ns{0}.
\end{align}

One can then define the coalgebra $C\rt D$, with $C\ot D$
 as underlying vector a space and the following coalgebra structure:
\begin{align}\label{cop}
&\D (c\rt d )= c\ps{1}\ot c\ps{2}\ns{-1}d\ps{1}\ot c\ps{2}\ns{0}\rt d\ps{2}\\\label{counit}
&\epsilon(c\rt d)=\epsilon(c)\epsilon(d).
\end{align}
It is straightforward to check that $C\rt D$ is a counital coassociative  coalgebra.

\smallskip

In addition we shall assume that
$C$ satisfies  the following compatibility condition, relating
the coaction of $\Kc$ and the action of $\Hc$:
\begin{equation}\label{compat}
\rc (hc)=h\ns{-1}c\ns{-1}\ot h\ns{0}c\ns{0}, \quad  \forall h\in \Hc, \, c\in C.
\end{equation}

One endows $C\rt D$ with a left action of $\Hc\rt \Kc$ by defining,
for any $h\in \Hc$, $k\in \Kc$, $c\in C$ and $d\in D$,
\begin{equation}\label{action}
(h\rt k)(c\rt d)=hc\rt kd.
\end{equation}
\begin{proposition}
The action defined by {\rm (\ref{action})} makes $C\rt D$ a $\Hc\rt\Kc$ module coalgebra.
\end{proposition}
\begin{proof} Using (\ref{cop}), (\ref{compat}), (\ref{c-c-1}), and the fact that $K$ is commutative, we have:
\begin{align*}
&\Delta ((h\rt k)(c\rt d))=\Delta(hc\rt kd)= \\
&(hc)\ps{1}\rt (hc)\ps{2}\ns{-1}(kd)\ps{1}\ot (hc)\ps{2}\ns{0}\ot (kd)\ps{2}=\\
&h\ps{1}c\ps{1}\rt h\ps{2}\ns{-1}c\ps{2}\ns{-1}k\ps{1}d\ps{1}\ot h\ps{2}\ns{0}c\ps{2}\ns{0}\ot k\ps{2}d\ps{2}=\\
&(h\rt k)\Delta(c\rt d).
\end{align*}
\end{proof}
Now let $M$ and $N$ be right-left SAYD modules  over $\Hc$  and $\Kc$ respectively.
\begin{definition}
We say that the SAYD module $M$ is $\Kc$-coinvariant if for any $h\in \Hc$, and any $m\in M$:
\begin{align}\label{coinv1}
& 1) \quad h\ns{-1}\ot mh\ns{0}=1\ot mh, \\ \label{coinv2}
& 2)\quad  m\sns{-1}\ns{-1}\ot m\sns{-1}\ns{0}\ot  m\sns{ 0} =1\ot m\sns{-1}\ot m\sns{0}
\end{align}
\end{definition}
\begin{definition}
We say that the SAYD module $N$ is $\Hc$-stable if for any  $h\in \Hc$, and $n\in N$ one has:
\begin{equation}\label{stable}
nh\ns{-1}\ot h\ns{0}=n\ot h
\end{equation}
\end{definition}

We now endow $M\ot N$ the a right action and left  coaction of $\Hc\rt \Kc$ as follows
\begin{align}\label{raction}
&(M\ot N)\ot \Hc\rt\Kc\ra (M\ot N), \\
&(m\ot n)(h\rt k)=mh\ot nk\\
&\notag \\ \label{lcoaction}
&(M\ot N)\ra \Hc\rt \Kc\ot (M\ot N), \\
&(m\ot n)\mapsto m\sns{-1}\rt \ot n\sns{-1}\ot  (m\sns{0}\ot n\sns{0})
\end{align}

\begin{proposition}
Let $M$ be $\Kc$-coinvariant  and $N$ be $\Hc$-stable.
Then via  the action and coaction defined by  {\rm  (\ref{raction})} and {\rm (\ref{lcoaction})} respectively, $M\ot N$ is an SAYD module over $\Hc\rt \Kc$.
\end{proposition}
\begin{proof}
First we show  that $M\ot N$ is a comodule over $\Hc\rt\Kc$.

 Indeed, using  (\ref{coinv2}) and the fact that $M$ and $N$ are comodule over $\Hc$, and $\Kc$ respectively, one has
\begin{align*}
& (m\ot n)\sns{-1}\ps{1}\ot (m\ot n)\sns{-1}\ps{2}\ot(m\ot n)\sns{0}=\\
&(m\sns{-1}\rt  n\sns{-1})\ps{1}\ot  (m\sns{-1}\rt \ot n\sns{-1})\ps{2} \ot m\sns{0}\ot n\sns{0}=\\
&m\sns{-1}\ps{1}\rt m\sns{-1}\ps{2}\ns{-1}n\sns{-1}\ps{1}\ot m\sns{-1}\ps{2}\ns{0}\rt n\sns{-1}\ps{2}\ot m\sns{0}\ot n\sns{0}=\\
&m\sns{-2}\rt m\sns{-1}\ns{-1}n\sns{-1}\ps{1}\ot m\sns{-1}\ns{0}\rt n\sns{-1}\ot m\sns{0}\ot n\sns{0}=\\
&m\sns{-2}\rt n\sns{-2}\ot m\sns{-1}\rt n\sns{-1}\ot m\sns{0}\ot n\sns{0}=\\
&(m\ot n)\sns{-2}\ot (m\ot n)\sns{-1}\ot (m\ot n)\sns{0}.\\
\end{align*}
Let us prove that $M\ot N$  satisfies (\ref{SAYD}).
First of all one sees that
\begin{align*}
&(\Delta\ot id)(\D(h\rt k))=\\
&h\ps{1}\rt h\ps{2}\ns{-1}h\ps{3}\ns{-2}k\ps{1}\ot
h\ps{2}\ns{0}\rt h\ps{3}\ns{-1}k\ps{2}   \ot  h\ps{3}\ns{0}\rt k\ps{3}.
\end{align*}
Using the above identity, the fact that both $M$ and $N$ satisfy (\ref{SAYD}),
that $\Kc$ is commutative,  that $M$ satisfies (\ref{stable}), and that $N$ satisfies (\ref{coinv1}),  
one can write successively 
\begin{align*}
&S((h\rt k)\ps{3})(m\ot n)\sns{-1}(h\rt k)\ps{1}\ot (m\ot n)\ns{0}(h\rt k)\ns{2}=\\
&S(h\ps{3}\ns{0}\rt k\ps{3})(m\sns{-1}\rt n\sns{-1})(h\ps{1}\rt h\ps{2}\ns{-1}h\ps{3}\ns{-2}k\ps{1})\ot \\
&\ot(m\sns{0}\ot n\sns{0})(h\ps{2}\ns{0}\rt h\ps{3}\ns{-1}k\ps{2})=\\
& S(h\ps{3}\ns{0})m\sns{-1}h\ps{1}\rt S(h\ps{3}\ns{-1}k\ps{3})n\sns{-1}h\ps{2}\ns{-1}h\ps{3}\ns{-3}k\ps{1}\ot\\
&\ot m\sns{0}h\ps{2}\ns{0}\ot n\sns{0}h\ps{3}\ns{-2}k\ps{2}=
\end{align*}
\begin{align*}
&S(h\ps{3}\ns{0})m\sns{-1}h\ps{1}\rt S(h\ps{3}\ns{-1}k\ps{3})n\sns{-1}h\ps{3}\ns{-2}k\ps{1}\ot\\
&\ot m\sns{0}h\ps{2}\ot n\sns{0}k\ps{2}=\\
& S(h\ps{3}\ns{0})m\sns{-1}h\ps{1}\rt S(h\ps{3}\ns{-1})h\ps{3}\ns{-2}S(k\ps{3})n\sns{-1}k\ps{1}\ot\\
&\ot m\sns{0}h\ps{2}\ot n\sns{0}k\ps{2}=\\
&S(h\ps{3})m\sns{-1}h\ps{1}\rt S(k\ps{3})n\sns{-1}k\ps{1}\ot m\sns{0}h\ps{2}\ot n\sns{0}k\ps{2}=\\
&(mh)\sns{-1}\rt (nk)\sns{-1}\ot (mh)\sns{0}\ot (nk)\sns{0}\\
&=[(m\ot n)(h\rt k)]\sns{-1}\ot [(m\ot n)(h\rt k)]\sns{0}.
\end{align*}
It is easily seen  that if $M$ and  $N$ are stable, then so is $M\ot N$.
\end{proof}

 For convenience, we recall (see \cite{hkrs2})  the definition of Hopf cyclic cohomology of a
left $\Hc$-module coalgebra with coefficients in a SAYD module $M$.
 It is the cyclic cohomology of  the cocyclic module defined by
 \begin{equation*}
 \Cc^n=C^n_H(C,M)=M\ot_H C^{\ot n+1}, \quad n\ge 0,
 \end{equation*}
 with the following cocyclic structure, where we abbreviate $\td c=c^0\odots c^n$,
 \begin{align}\label{CCM1}
& \p_i(m\ot_\Hc \td c )=m\ot_\Hc c^0\odots \Delta(c_i)\odots c^n, \\\label{CCM2}
 &\p_{n+1}(m\ot_\Hc \td c)=m\sns{0}\ot_\Hc c^0\ps{2}\ot c^1\odots c^n\ot m\sns{-1}c^0\ps{1},\\\label{CCM3}
 &\sigma_i(m\ot_\Hc \td c)=m\ot_\Hc c^0\odots \epsilon(c^{i+1})\odots c^n, \\\label{CCM4}
 &\tau(m\ot_\Hc \td c)=m\sns{0}\ot_\Hc c^1\odots c^n\ot m\sns{-1}c^0.
 \end{align}

 We define the bigraded module $\FX^{p,q}$ by
 \begin{equation*}
 \FX^{p,q}= N\ot_\Kc D^{\ot p+1}\ot   M\ot_{\Hc}\ot C^{\ot q+1} \, ,
  \end{equation*}
 and endow $\FX$ wit the operators
  \begin{align}
&\hd_i:\FX^{(p,q)}\rightarrow \FX^{(p+1,q)}, && 0\le i\le p+1\\
&\hs_j: \FX^{(p,q)}\rightarrow \FX^{(p-1,q)}, &&0\le j\le p-1\\
&\hta:  \FX^{(p,q)}\rightarrow \FX^{(p,q)},
\end{align}
defined by
\begin{align}\label{chdi}
&\hd_j(n\ot\td{d}\ot   m\ot \td{c})= n\ot c^0\ot\dots\ot \Delta(c^i)\ot\dots \ot c^p\ot   m\ot \td{d},\\ \label{chdn}
&\hd_{p+1}(n\ot\td{d}\ot   m\ot \td{c})=\\ \notag
&n\sns{0}\ot d^0\ps{2}\odots  d^p\ot \td{c}\ns{-1}n\sns{-1}d^0\ps{1} \ot   m\ot \td{c}\ns{0},\\ \label{chs}
&\hs_j(n\ot\td{d}\ot   m\ot \td{c})=n\ot d^0\odots \epsilon(d^{j+1})\odots d^p\ot   \td{c},\\ \label{cht}
&\hta(n\ot\td{d}\ot   m\ot \td{c})=\\ \notag
& n\sns{0}\ot d^1\odots  d^p\ot \td{c}\ns{-1}n\sns{-1}d^0 \ot   m\ot \td{c}\ns{0}  ;
\end{align}
the vertical structure is just the cocyclic structure of $\Cc^\ast$, with
\begin{align} \label{cvd}
&\vd_i= \Id\ot   \p_i:\FX^{(p,q)}\rightarrow \FX^{(p,q+1)}, &&0\le i\le q+1\\ \label{cvs}
&\vs_j=\Id\ot   \sigma_j: \FX^{(p,q)}\rightarrow \FX^{(p,q-1)}, && 0\le j\le q-1\\ \label{cvt}
&\vta= \Id\ot   \tau:  \FX^{(p,q)}\rightarrow \FX^{(p,q)},
\end{align}
where $\p_i$, $\sigma_i$, and $\tau$ are  defined by (\ref{CCM1})-(\ref{CCM4}).

\begin{proposition}
The bigraded module $\FX$ endowed with the above operators
$(\hd,\hs,\hta, \vd,\vs,\vta)$ defines a bicocyclic module.
\end{proposition}
 \begin{proof}
 The $p$th column of  $\FX$ is identical to
 $N\ot_\Kc D^{\ot p}\ot   C^\ast_\Hc(C,M)$ and the vertical operators do not
 affect the horizontal part. So the columns are cocyclic modules.

 To see that the rows are also cocyclic modules,
 one identifies the $q$th row  with $C^ \ast_{\Kc}(D,N\ot M\ot C^{\ot q+1})$,
  where $\Kc$ acts on  $N\ot M\ot C^{\ot q+1}$ via its action on $N$, and coacts as follows.
  For  $n\ot m\ot \td c\in N\ot M\ot C^{\ot q+1}$ the coaction is
 \begin{equation*}
 n\ot m\ot \td c\mapsto \td c\ns{-1} m\sns{-1}\ot n\sns{0}\ot m\ot \td c\ns{0}.
 \end{equation*}
 Using the fact that $N$ is a SAYD module over $\Kc$ and that $\Kc$ is commutative,
 one sees that $N\ot M\ot C^{\ot q+1} $ is a right-left SAYD module over $\Kc$, hence
 the row forms a cocyclic module.

 To finish the proof it remains to show that the vertical and horizontal operators
 commute among each other. We shall only check
 the nontrivial part, that is the commutation
 of $\vta$ and $\hta$.
 Using (\ref{coinv2}) one can write
 \begin{align*}
 &\hta\vta(n\ot \td d \ot   m\ot \td c)=\\
 &\hta( n\ot \td d\ot   m\sns{0}\ot c^1\odots c^q \ot m\sns{-1} c^0 )=\\
 &n\sns{0}\ot d^1\odots d^p\ot  n\sns{-1}m\sns{-1}\ns{-1}\td c \ns{-1} d^0\ot   \\
 & m\sns{0}\ot c^1\ns{0}\odots c^q\ot m\sns{-1}\ns{0}c^0=\\
 &n\sns{0}\ot d^1\odots  d^p\ot \td{c}\ns{-1}n\sns{-1}d^0 \ot   \\
  &m\sns{-1}\ot c^1\ns{0}\odots c^q\ns{0}m\sns{-1}c^0\ns{0}=\\
   &\vta(n\sns{0}\ot d^1\odots  d^p\ot \td{c}\ns{-1}n\sns{-1}d^0 \ot   m\ot \td{c}\ns{0} )=\\
 &\vta\hta(n\ot \td d \ot   m\ot \td c).
 \end{align*}

 \end{proof}
We denote the diagonal of   $\FX$ by $\FY$;
it is a cocyclic module whose degree $n$ component is
$$N\ot_\Kc\ot D^{\ot n}\ot   M\ot_\Hc C^{\ot n},$$
with the following operators:
\begin{align}
&d_i=\hd_i \vd_i,&  0\le i \le n+1,\\
& s_j=\hs_j\vs_j,& 0\le j\le n-1,\\
&t=\hta\vta
\end{align}
  Next we identify the cocyclic module $C^\ast_{\Hc\rt \Kc}(C\rt D, M\ot N)$ with $\FY$.

 \begin{proposition}
 The following map
 \begin{equation*}
\Psi:   C^\ast_{\Hc\rt \Kc}(C\rt D, M\ot N)\ra  \FY
\end{equation*}
 is well-defined and defines an isomorphism of cocyclic modules:
 \begin{align*}
\Psi&(m\ot n\ot c^0\rt d^0\ot \dots\ot c^n\rt d^n)=
 n\ot d^0 \ot c^1\ns{-n}d^1\ot\dots\\
&\ot c^1\ns{-2}\dots c^{n-1}\ns{-2}d^{n-1}
\ot  c^1\ns{-1}\dots c^n\ns{-1}d^n\\
&\ot   m\ot c^0\ot c^1\ns{0}\odots
c^n\ns{0}.
\end{align*}
\end{proposition}
\smallskip

\begin{proof}
One can check that the following expression defines an inverse for
$\Psi$:
\begin{align*}
&\Psi^{-1}(n\ot\tilde{d}\ot   m\ot \tilde{c})=m\ot n\ot c^0\rt d^0 \ot c^1\ns{0}\rt S(c^1\ns{-1})d^1\ot\dots \\
&\dots \ot c^{n}\ns{0}\rt S(c^1\ns{-n}\dots
c^n\ns{-1})d^n.
\end{align*}
First we show that $\Psi$ is well-defined. Let  $h\in \Hc$, $k\in \Kc$,  $\wt{h\rt k}\in (\Hc\rt\Kc)^{\ot n+1}$, $m\in M$ and $n\in N$.
One observes that
\begin{align*}
&\D^{(n)}(h\rt k)= \\
&h\ps{1}\rt h\ps{2}\ns{-1}\dots h\ps{n+1}\ns{-n}k\ps{1}\ot h\ps{2}\ns{0}\rt h\ps{3}\ns{-1}\dots h\ps{n+1}\ns{-n+1}k\ps{2}\ot\dots\\
&\dots \ot  h\ps{n}\ns{0}\rt h \ns{n+1}\ns{-1}k\ps{n}\ot h\ps{n+1}\ns{0}\rt k\ps{n+1}.
\end{align*}
Then by using (\ref{stable}) we have:
\begin{align*}
&\Psi(m\ot n\ot (h\rt k)(\wt{c\rt d}))= \Psi (m\ot n\ot \D^{(n)}(h\rt k)(\wt{ c\rt d}) =\\
&\Psi(m\ot n\ot h\ps{1}c^0\rt h\ps{2}\ns{-1}\dots h\ps{n+1}\ns{-n}k\ps{1}d^0\ot h\ps{2}\ns{0}c^1\rt h\ps{3}\ns{-1}\dots\\
& h\ps{n+1}\ns{-n+1}k\ps{2}d^1\odots  h\ps{n+1}\ns{0}c^n\rt k\ps{n+1}d^n)=\\
&n\ot   h\ps{2}\ns{-n}\dots h\ps{n+1}\ns{-n}k\ps{1} c^1\ns{-n}d^0 \ot\dots \\
&\dots \ot  h\ps{2}\ns{-n}\dots h\ps{n+1}\ns{-n}k\ps{n+1} c^1\ns{-1}\dots c^n\ns{-1} d^n\ot  \\
& m\ot  h\ps{1}c^0\ot h\ps{2}\ns{0} c^1\ns{0}\odots h\ps{n+1}\ns{0}c^n\ns{0}=\\
&nh\ps{2}\ns{-1}\dots h\ps{n+1}\ns{-1}\ot k\ps{1}c^1\ns{-n}d^0 \odots   k\ps{n+1}c^1\ns{-1}\dots c^n\ns{-1} d^n\ot  \\
& m\ot  h\ps{1}c^0\ot h\ps{2}\ns{0} c^1\ns{0}\odots h\ps{n+1}\ns{0}c^n\ns{0}=\\
&n\ot  k\ps{1}c^1\ns{-n}d^0 \odots   k\ps{n+1}c^1\ns{-1}\dots c^n\ns{-1} d^n\ot   \\
&m\ot  h\ps{1}c^0\ot h\ps{2} c^1\ns{0}\odots h\ps{n+1}c^n\ns{0}=\\
&nk\ot  c^1\ns{-n}d^0 \odots   c^1\ns{-1}\dots c^n\ns{-1} d^n\ot  \\
& mh\ot  c^0\ot h\ps{2} c^1\ns{0}\odots c^n\ns{0}=  \Psi(mh\ot nk\ot \wt{c\rt d}). \\
\end{align*}
The next task
is to show that $\Psi$ is a cyclic map. To do this we should check that $\Psi$ commutes with
the faces, the degeneracies and the cyclic operator.
We will just check the commutativity of $\Psi$ with $\tau$ and $\p_0$,
the other verifications being completely similar.

Indeed, by using (\ref{stable}), and (\ref{coinv2}) one can write
\begin{align*}
&\Psi(\tau(m\ot n\ot \wt{c\rt d}))= \\
&\Psi (m\sns{0}\ot n\sns{0}\ot c^1\rt d^1\odots c^n\rt d^n\ot m\sns{-1}c^0\rt n\sns{-1}d^0)=\\
& n\sns{0}\ot d^1\ot c^2\ns{-n}d^2\odots c^2\ns{-1}\dots c^n\ns{-1}m\sns{-1}\ns{-1}c^0\ns{-1}d^0\ot  \\
& m\sns{0}\ns{0}\ot c^1\ot c^2\ns{0}\odots c^n\ns{0}\ot m\sns{-1}\ns{0}c^0\ns{0}=\\
&n\sns{0}\ot d^1\ot c^2\ns{-n}d^2\odots c^2\ns{-1}\dots c^n\ns{-1}c^0\ns{-1}d^0\ot  \\
& m\sns{0}\ot c^1\ot c^2\ns{0}\odots c^n\ns{0}\ot m\sns{-1}c^0\ns{0}=\\
&n\sns{0}c^1\ns{-1}\ot d^1\ot c^2\ns{-n}d^2\odots c^2\ns{-1}\dots c^n\ns{-1}c^0\ns{-1}d^0\ot  \\
& m\sns{0}\ot c^1\ns{0}\ot c^2\ns{0}\odots c^n\ns{0}\ot m\sns{-1}c^0\ns{0}=\\
&n\sns{0}\ot c^1\ns{-n-1} d^1\ot c^1\ns{-n}c^2\ns{-n}d^2\odots c^1\ns{-1}c^2\ns{-1}\dots c^n\ns{-1}c^0\ns{-1}d^0\ot  \\
& m\sns{0}\ot c^1\ns{0}\ot c^2\ns{0}\odots c^n\ns{0}\ot m\sns{-1}c^0\ns{0}=\\
& t \Psi(m\ot n\ot \wt{c\rt d})).
\end{align*}
Fnally, using once more  (\ref{stable}) we check that $\Psi\p_0=d_0\Psi$,
as follows:
\begin{align*}
&\Psi\p_0(m\ot n\ot \wt{c\rt d})=\\
& \Psi(m\ot n\ot c^0\ps{1}\rt c^0\ps{2}\ns{-1}d^0\ps{1}\ot c^0\ps{2}\ns{0}\rt d^0\ps{2}\ot c^1\rt d^1\odots c^n\rt d^n)=\\
&n\ot c^0\ps{2}\ns{-n-1}d^0\ps{1}\ot c^0\ps{2} \ns{-n}d^0\ps{2}\ot c^0\ps{2}\ns{-n+1}c^1\ns{-n+1}d^1\ot\dots\\
&\dots \ot  c^0\ps{2}\ns{-1} c^1\ns{-1}\dots c^n\ns{-1}d^n\ot    c^0\ps{1}\ot c^0\ps{2}\ns{0}\odots c^n\ns{0}=\\
&n c^0\ps{2}\ns{-1}\ot  d^0\ps{1}\ot d^0\ps{2}\ot c^1\ns{-n+1}d^1\ot\dots\\
&\dots \ot  c^1\ns{-1}\dots c^n\ns{-1}d^n\ot    c^0\ps{1}\ot c^0\ps{2}\ns{0}\odots c^n\ns{0}=\\
&n\ot  d^0\ps{1}\ot d^0\ps{2}\ot c^1\ns{-n+1}d^1\ot\dots\\
&\dots \ot  c^1\ns{-1}\dots c^n\ns{-1}d^n\ot    c^0\ps{1}\ot c^0\ps{2}\ot c^3\ns{0}\odots c^n\ns{0}=\\
& d_0\Psi(m\ot n\ot \wt{c\rt d)}.\\
 \end{align*}
\end{proof}

We have thus proved the following generalization of Theorem \ref{mix}.

  \begin{theorem} \label{mix+}
  The mixed complexes  $(C^*_{\Hc\rt \Kc}(C\rt D,M\ot N), \,b,\, B)$  and $(\Tot(\FX), \,b_T,\, B_T))$ are quasi-isomorphic.
  \end{theorem}


 \section{A Cartan homotopy formula for Hopf cyclic
 cohomology}\label{S3}

In this section we prove a homotopy  formula for Hopf cyclic
cohomology of coalgebras with coefficients in SAYD modules, which
in conjunction with the spectral sequence of the previous section
will allow us to compute the Hopf cyclic cohomology of Hopf algebras
such as $\Hc_1^\dag$ and $\Hc_1$. To put it into the proper perspective,
we begin by recalling the cyclic cohomological
version of the Cartan-type homotopy formula (cf.~\cite[Part II, Prop. 5]{cng},
\cite[\S 4.1]{lo})
for actions of derivations on algebras.
\medskip

 Let $\Ac$ be an algebra and let
 $D: \Ac \rightarrow \Ac$ be a derivation. Then $D$ extends to a
 `Lie derivative' operator,
 acting on the cochain complex of $\Ac$ by
$$\Lc_Df(a^0\ot\dots \ot a^n)=
\sum_{i=0}^{n}f(a_0\ot \dots \ot a_{i-1}\ot D(a_i)\ot a_{i+1} \dots \ot a_n) .
$$
The `contraction' by $D$ is made from the operators
 $$e_D: C^n(\Ac)\rightarrow C^{n+1}(\Ac) \quad \text{and} \quad
 E_D:
C^n(\Ac)\rightarrow C^{n-1}(\Ac)
$$
defined by
\begin{align*}
&e_Df(a_0\ot\dots \ot a_{n+1})= (-1)^{n}f(D(a^n)a^0\ot\dots \ot a^{n-1}), \\\notag
&\\
&E_Df(a_0\ot  \dots\ot a_{n-1})=\\
& \sum_{1\le i\le j\le n}(-1)^{n(i+1)}f(1\ot a_i\ot a_{i+1}\ot
\dots\ot a_{j-1}\ot D(a_j)\ot a_{j+1}\ot \dots\\ &\dots \ot a_n\ot a_0\ot
\dots\ot a_{i-1} )
\end{align*}
and the following relations are satisfied (cf.~\cite[\S 4.1]{lo}):
\begin{equation*}
[e_D,b]=0\, , \qquad [E_D,B]=0 \, ,
\end{equation*}
\begin{equation*}
[e_D,B]+[E_D,b]=\Lc_D \, .
\end{equation*}
\medskip

As an example, we can take $\Ac=C^{\infty}_c(F^+M)\rtimes\Gamma$ as in
$\S$\ref{S1.1}, with $D=Y\in \Hc_1$ acting as a
derivation on $\Ac$ .

We can pullback the above operators on $C^n(\Hc_1)$ via the characteristic map
$$\chi: \Hc_1^{\ot n}\rightarrow Hom(\Ac^{\ot n+1},\bbc),$$
$$\chi(h^1\odots  h^n)(a_0\odots\a_n)=\tau(a_0h^1(a_1)\dots h^n(a_n))
$$
and obtain operators $\Lc_Y^{tr}$, $e_Y^{tr}$ and
$E_Y^{tr}$ respectively. The first one can be intrinsically expressed
as follows.

\begin{lemma}
In the above notation,
 $$\Lc_Y^{tr}=\Id-\ad Y: \Hc^{\ot n}\rightarrow
\Hc^{\ot n }.
$$
\end{lemma}

\begin{proof}
Recalling that the invariant trace $\tau$ satisfies
$\tau(h(a)b)=\tau(aS_{\delta}(h)b)$, and $S_\delta(Y)=-Y+1$, one obtains
 \begin{align*}
&\Lc_Y\chi(h^1\odots h^n)(a_0\odots a_n)=\\
&=\chi(h^1\odots h^n)(Y(a_0)\ot a_1\odots a_{n})+ \\
&+\sum_{i>0}\chi(h^1\odots h^n)(a_0\odots \ot a_{i-1}\ot Y(a_i)\ot a_{i+1}\odots a_n)\\
&=\tau(Y(a_0)h^1(a_1)\dots h^n(a_n))+\sum_{i>0}\tau(a_0h^1(a_1)\dots h^i(Y(a_i)) \dots h^n(a_n))\\
&=\tau(a_0S_\delta(Y)(h^1(a_1)\dots h^n(a_n)) +\sum_{i>0}\tau(a_0h^1(a_1)\dots h^i(Y(a_i)) \dots h^n(a_n))\\
&=\tau(a_0h^1(a_1)\dots h^n(a_n))-\sum_{i>0}\tau(a_0h^1(a_1)\dots Y(h^i(a_i))\dots h^n(a_n) )+\\
&+\sum_{i>0}\tau(a_0h^1(a_1)\dots h^i(Y(a_i)) \dots h^n(a_n))\\
&=\tau(a_0h^1(a_1)\dots h^n(a_n))-\sum_{i>0}\tau(a_0h^1(a_1)\dots
[Y,h^i](a_i))\dots h^n(a_n)).
\end{align*}
\end{proof}
\smallskip

Similarly, one can obtain intrinsic expressions for the pullbacks of the operators
$e_Y^{tr}$ and $E_Y^{tr}$.

Guided by this example, we proceed to find the
analogous homotopy formula in the general context of
 Hopf cyclic cohomology of coalgebras with coefficients in  SAYD modules.
\medskip

Let  $D: C\rightarrow C$  be an $\Hc$-linear
 coderivation, that is an $\Hc$-linear map such that
 $$\Delta(D(c))=
 D(c\ps {1})\ot c\ps{2}+c\ps{1}\ot D(c\ps{2}).
 $$
 We define the
 operators   $\Lc_D: \Cc^n\ra \Cc^n$, $e_D:\Cc^n\ra \Cc^{n+1}$ and
 $E_D: \Cc^n\ra \Cc^{n-1}$ on
 $\, \Cc^n$ as follows:
\begin{align*}
&\Lc_D(m\ot c^0\odots c^n)=\sum_{i=0}^{n}m\ot c^0\odots D(c_i)\odots
c^n, \\
&e_D(m\ot c^0\odots c^n)=(-1)^{n}m\sns{0}\ot c^0\ps{2}\ot c^1\odots
c^n\ot m\sns{-1}D(c^0\ps{1}), \\
&E_D=\sum_{j=1}^{n}\sum_{i=1}^{j} E^{j,i}_D: \Cc^n \ra \Cc^{n-1}\, ,
\end{align*}
where
\begin{align*}
&E_D^{j,i}(m\ot c^0\odots c^{n})=\\
&=(-1)^{n(i+1)}\epsilon(c^0)m\sns{0}\ot c^{n-i+2}\odots c^{n+1}\ot
m\sns{-1}c_1\odots m\sns{-(j-i)}c^{j-i}\ot
\\
&\ot m\sns{-(j-i+1)}D(c^{j-i+1})\ot m\sns{-(j-i+2)}c^{j-i+2}\odots
m\sns{-(n-i+1)}c^{n-i+1} .
\end{align*}
Let us also introduce the auxiliary  operators $\, \psi_j:\Cc^n\ra \Cc^n$,
 \begin{equation*}
 \psi_j(m\ot c^0\odots c^n)= m\ot c^0\odots D(c_j)\odots c^n\, .
\end{equation*}

\begin{lemma}
 The operators $\Lc_D$,
$e_D$ and $E_D$ are well-defined and the following identities hold:
\begin{align*}
&\Lc_D=\sum_{j=0}^n\psi_j,\\
&e_D=(-1)^{n}\tau\psi_0\p_0=(-1)^{n}\psi_{n+1}\p_{n+1}\,\\
&E_D=\sum_{j=1}^{n}\sum_{i=1}^{j}(-1)^{n(i+1)}\psi_j\tau^{-i}\sigma_n\tau.
\end{align*}
\end{lemma}
\begin{proof}
The operators $\Lc_D$, $E_D$, and $e_D$ are well-defined as soon as
 $\psi_j$s are shown to be well-defined. The latter property is easily
checked, because $D$ is $\Hc$-linear.

The first two identities being
obviously  true, we only check the last one:
\begin{align*}
&\psi_j\tau^{-i}\sigma_n\tau(m\ot c^0\odots c^n)=
\psi_j\tau^{-i}(\epsilon(c_0)m\ot c^1\odots c^n)\\
&=\psi_j(\epsilon(c^0) m\sns{0}\ot S^{-1}(m\sns{-1}) c^{n-i+2}\odots
S^{-1}(m\sns{-i})c^{n+1}\ot c_1\ot\dots\\
&\dots\ot c^{n-i+1}\\
&= \epsilon(c^0)m\sns{0}\ot c^{n-i+2}\odots c^{n+1}\ot
m\sns{-1}c_1\odots m\sns{-(j-i)}c^{j-i}\ot
\\
&\ot m\sns{-(j-i+1)}D(c^{j-i+1})\ot m\sns{-(j-i+2)}c^{j-i+2}\odots
m\sns{-(n-i+1)}c^{n-i+1}.\\
\end{align*}
\end{proof}
\begin{lemma}\label{e-E}
The operator $e_D$  and $E_D$  commute, in the graded sense,
with the Hochschild coboundary  $b$ and
with the Connes boundary  $B$ respectively:
\begin{equation*}
[b,e_D]=0, \qquad  [B, E_D]=0.
\end{equation*}
\end{lemma}
\begin{proof}
The commutation of $E_D$ and $B$ holds
because of the fact that we are always working in the normalized subcomplex.

For $e_D$ and $b$, one sees that $\p_i e_D=e_D\p_i$,  for $0\le i\le n$, then by
using the coderivation property of $D$ one checks that $\p_{n+1}
e_D-\p_{n+2}e_D=e_D\p_{n+1}$.
\end{proof}

\begin{lemma}\label{E_D}
 $$e_DB=-\sum_{k=1}^{n}(-1)^{(n+1)(n-k+2)}(-1)^{k}E_D^{n,n-k+1}\p_k.$$
\end{lemma}
\begin{proof}
\begin{align*}
&e_DB=\sum_{i=0}^{n-1}(-1)^{n-1}(-1)^{ni}\psi_n\p_n\tau^i\sigma_{n-1}\tau
=\sum_{i=0}^{n-1}(-1)^{ni+n-1}\psi_n\tau\p_0\tau^i\sigma_{n-1}\tau\\
&=\sum_{i=0}^{n-1}(-1)^{ni+n-1}\psi_n\tau^i\p_i\sigma_{n-1}\tau
=\sum_{i=0}^{n-1}(-1)^{ni+n-1}\psi_n\tau^{i+1}\sigma_n\p_i\tau\\
&=\sum_{i=0}^{n-1}(-1)^{ni+n-1}\psi_n\tau^{i+1}\sigma_n\tau\p_{i+1}
=-\sum_{k=1}^{n}(-1)^{(n+1)(n-k+2})(-1)^{k}\psi_n\tau^{k}\sigma_n\tau\p_{k}\\
&=-\sum_{k=1}^n (-1)^{(n+1)(n-k+2)}(-1)^{k}E_D^{n,n-k+1}\p_k.\\
\end{align*}
\end{proof}

\begin{lemma}\label{e_D}
$$Be_D=\sum_{i=0}^{n}(-1)^{n+ni}\psi_{n-i}\tau^{i+1}.$$
\end{lemma}
\begin{proof}
\begin{align*}
&Be_D=\sum_{i=0}^{n}(-1)^{n+ni}\tau^i\sigma_n\tau\psi_{n+1}\p_{n+1}
=\sum_{i=0}^{n}(-1)^{n+ni}\tau^i\sigma_n\psi_{n}\tau\p_{n+1}\\
&=\sum_{i=0}^{n}(-1)^{n+ni}\tau^i\psi_{n}\sigma_n\p_{n}\tau
=\sum_{i=0}^{n}(-1)^{n+ni}\psi_{n-i}\tau^{i+1}
\end{align*}
\end{proof}
\begin{lemma}\label{L_D}
$$\Lc_D=Be_D+E_D\p_0+(-1)^{n+1}E_D\p_{n+1}.$$
\end{lemma}
\begin{proof}
Using Lemma \ref{e_D}, one obtains
\begin{align*}
&Be_D+E_D\p_0+(-1)^{n+1}E_D\p_{n+1}= \sum_{i=0}^{n}(-1)^{n+ni}
\psi_{n-i}\tau^{i+1}+\\
&+\sum_{j=1}^n \sum_{i=1}^j
(-1)^{(n+1)}E_D^{j,i}\p_{n+1}+E_D\p_{0}=(\psi_0+\sum_{j=1}^{n}E_D^{j,1}\p_{n+1})+\\
&+(\sum_{i=0}^{n-1} (-1)^{n+ni}\psi_{n-i}\tau^{i+1}+ \sum_{j=1}^n \sum_{i=2}^j
 (-1)^{n+1}E_D^{j,i}\p_{n+1}+E_D\p_{0}).
\end{align*}

\ni To finish the proof we just need to show that the result of the
first parentheses is $\Lc_D$ and that of the second one vanishes. To
this end, one checks  that
$$E_D^{j,i}\p_0=(-1)^{(n+1)i}\psi_{j}\tau^{n+1-i}\,,
$$
and
$$E_D^{j,i}\p_{n+1}=(-1)^{(n+1)(i+1)}\psi_{j}\tau^{n+2-i}.
$$
\end{proof}

\begin{proposition}\label{homotop}
$$[e_D+E_D,b+B]=\Lc_D.$$
\end{proposition}
\begin{proof}
From Lemma \ref{e-E} we already know that $[E_D,B]=0$ and $[e_D,b]=0$.
To complete the proof,  after using Lemmas   \ref{E_D}, \ref{e_D}, and \ref{L_D}, it
suffices to show that:
$$ \p_kE^{j,i}=  \left\{
                       \begin{array}{ll}
                      E_D^{j,i}\p_{k-i+1}   , & { k=i+1,\dots,j-1;} \\
                        E^{j,i}_D\p_{k-i+1} , & {k=j+1\dots n;} \\
                         E^{j,i}\p_{k-i+1}+E^{j+1,i}\p_{k-i+1}, & {k=j;} \\
                         E^{j+1,i}\p_{n-i+k}, & {k=1,\dots,i-1;} \\
                         E^{j+1,i}\p_1, & \hbox{k=i.}
                       \end{array}
                     \right.
 $$
We check the first case and leave the rest to the reader. Indeed,
\begin{align*}
&E^{j,i_D}\p_{k-i+1}=
\psi_j\tau^{-i}\sigma_n\tau\p_{k-i+1}=\psi_j\tau^{-i}\sigma_n\p_{k-i}\tau=
\psi_j\p_k\tau^i\sigma_{n-1}\tau\\&=\p_k\psi_j\tau^{-i}\sigma_{n-1}\tau
=\p_kE^{j,i}_D.
\end{align*}
\end{proof}

Consider now a Hopf algebra $\Hc$  endowed with an MPI $(\delta,\sigma)$,
and let $Z \in \Hc $ be a primitive element, i.e.,
$$\Delta(Z) \, = \, Z\ot 1+1\ot
Z  \qquad \text{and} \qquad \epsilon(Z)\, =\, 0\, .
$$
Then $\, D_Z: \Hc \ra \Hc$ defined by
\begin{equation*}
D_Z(h)=hZ
\end{equation*}
 is an $\Hc$-linear coderivation, and
 we denote by $\Lc_Z$ its functorial extension to the tensor algebra
 of $\Hc$.

\begin{lemma}\label{ad}
$$\Lc_Z=\delta(Z) \Id- \ad Z$$
\end{lemma}
\begin{proof}
Let $\Theta: C^n_H(\Hc,~^\sigma\bbc_\delta)\ra  \Hc^{\ot n}$ be the
canonical isomorphism  defined by the assignment
\begin{equation*}
h^0\odots h^n\mapsto S_\delta(h^0)\cdot h^1\odots h^n \, ,
\end{equation*}
with inverse
$$h^1\odots h^n\mapsto 1\ot h^1\odots h^n.$$
Since $S_\delta(Z)=\delta(Z)-Z$,  one gets
\begin{align*}
&\Theta \Lc_Z \Theta^{-1}(h^1\odots h^n)=\Theta \Lc_Z(1\ot h^1\odots
h^n)=\\
& =\Theta(Z\ot h^1\odots h^n )+\sum_{i+1}^{n} \Theta(1\ot h^1\odots h^iZ\odots h^n )\\
&=\delta(Z)h^1\odots h^n\ot -\sum_{i+1}^n(h^1\odots Zh^i\odots
h^n)+\\
&+\sum_{i=1}^n h^1\odots h^iZ\odots h^n
=(\delta(Z)\Id-\ad_Z)(h^1\odots h^n).
\end{align*}
\end{proof}

To illustrate the application of this homotopy we consider the case
of $ \Hc_1$, with the primitive element
$Y\in \Hc_1$ in the role of $Z$.  First, from
Lemma \ref{ad} and
Proposition \ref{homotop}, one obtains:
\medskip

\begin{corollary}{\rm (The Cartan Homotopy Formula for $\Hc_1$)}\label{cartan}
$$ \Id-\ad Y\, =\, [E_Y+e_Y, B+b] \, .
$$
\end{corollary}
\medskip

We now look at the concrete way in which
 $\ad Y$ acts on $\Hc_1$. Since
$$\ad Y(Y)=0\, , \quad \ad Y(X)=X \quad \text{and} \quad
\ad Y(\delta_k)=k\delta_k\, , \quad \fl \, k \geq 1 \, ,
$$
one sees that, for a typical basis element
$Y^pX^q\delta_1^{r_1}\dots \delta_m^{r_m}$ of $\Hc_1$,  
$$\ad Y(Y^pX^q\delta_1^{r_1}\dots \delta_m^{r_m})=
(q+r_1\dots+r_m)Y^pX^q\delta_1^{r_1}\dots \delta_m^{r_m}.
$$
Thus,  $\ad Y$ is diagonalizable, acting as a grading operator, with non-negative
integer weights.  In particular, the weight $1$  eigenspace is
$$\Hc_1 [1] \, = \,
\bbc[Y]X\ot \bbc[Y]\delta_1.
$$
The above grading extends to the cocyclic module
${\Hc_1}_\nr$ and gives a corresponding decomposition
$${\Hc_1}_\nr =\bigoplus_{k\ge 0}{{\Hc_1}_\nr}[k], $$
where
$\displaystyle {{\Hc_1}_\nr}[k]= \bigoplus_{n \ge 0}\{\td{h}\in \Hc_1^{\ot n}\mid \nm{\td h}=k
\}$, with
$ |h^1 \ot \ldots \ot h^n| \, = \, |h^1| + \ldots +  |h^n| $
denoting the total weight.

Since the operators defining the
cocyclic structure of ${{\Hc_1}_\nr}$ obviously preserve the weight,
 the above decomposition is actually a decomposition of cocyclic modules.
 Then Corollary \ref{cartan} has in turn the following implications.

\begin{corollary} \label{percartan}
\begin{equation*}
HP^i_{(\delta,1)}({{\Hc_1}_\nr}[1]) = HP^i(\Hc_1) \quad
 {\rm and} \quad
HP^i_{(\delta,1)}({{\Hc_1}_\nr}[k])=0, ~~~ \fl \, k\neq 1.
\end{equation*}
\end{corollary}
\bigskip

\section{Hopf cyclic cohomology of Hopf algebras in codimension 1} \label{S4}

\subsection{Hopf cyclic cohomology of $\Hc_1$ and  $\Hc_{1 s}$}\label{S4.2}

By its very construction (cf. \cite{cm2}), the Hopf algebra
 $\Hc_1$ is
a bicrossed product of two Hopf subalgebras. We shall use this feature
to construct a spectral sequence converging to
the Hochschild  cohomology of $\Hc_1$. This spectral sequence is
further graded by the weight, and its  weight $1$ component is extremely
simple. This will allow us, first to explicitly compute the
weight $1$ component of the Hochschild  cohomology and then, in
conjunction
with Corollary \ref{percartan},  to finish off the calculation
of the periodic Hopf cyclic cohomology of $\Hc_1$.
 \medskip

To begin with, we shall make the bicrossed product construction
completely explicit (cf.~\cite{cm2}, also~\cite{kr}).
Let $\Uc=\Uc({\Fg}_{-})$  denote the universal enveloping algebra of
the Lie algebra ${\Fg}_{-}$ of the `$ax+b$'  group of affine
diffeomorphisms of the real line. The Lie algebra ${\Fg}_{-}$ has generators $X$, $Y$
with the bracket
\begin{equation}\label{XY}
[Y,X]=X.
 \end{equation}
Let $\Fc=\Fc(G_+)$ be the (commutative) Hopf subalgebra of  $\Hc_1$ generated
by $\{\delta_k \,  ; \, k \in \mathbb{N} \}$, which we may view as algebra of
regular functions on the pronilpotent Lie group $G_+$ consisting of the
diffeomorphisms
$\psi \in \Diff(\mathbb{R})$ that satisfy
$\psi(0)=0$ and $\psi'(0)=1$.

 The Hopf algebra $\Uc$ acts from the right on the Hopf algebra
 $\Fc$ via
 \begin{equation*}
  \delta_k\lt X=-\delta_{k+1},    \quad\delta_k\lt Y=-k\delta_k,
  \end{equation*}
  making $\Fc$ a  right $\Uc$ -module algebra.
  \smallskip

  On the other hand, the Hopf algebra $\Fc$ coacts from the left on $\Uc$ via
  the coaction $\rho: \Uc \rightarrow \Fc \ot \Uc$ defined by
  \begin{align}
  &\rho(X)=1\ot X+\delta_1\ot Y, \\
  & \rho(Y)=1\ot Y, \\
  &\rho(1)=1\ot 1,
  \end{align}
   and extends to the basis of $\Uc$ via
   \begin{equation}\label{bi}
  \rho(u^1u^2)= (u^1\ns{-1}\lt u^2\ps{1})u^2\ps{2}\ns{-1}\ot u^1\ns{0}u^2\ps{2}\ns{0} .
\end{equation}
With respect to this coaction $\Uc$ is a left $\Fc$-comodule coalgebra,
Moreover, it is not difficult to check that the above
 action and coaction satisfy the axioms \cite[page 232]{maj} required
 to form the bicrossed product. The resulting Hopf algebra coincides
 as coalgebra with the cocrossed product coalgebra  $\Uc \ltimes \Fc$,
 and as algebra with the crossed product  $ \Uc \rtimes \Fc$.
 We identify $X$ ,  $Y$ and $\delta_k \in \Hc_1$
 with $X\rt 1$, $Y\rt 1$ and $1\rt \delta_k$ in $\Uc \bi \Fc$, respectively.

 \begin{lemma} \label{bicp}
  $\Hc_1$ and $ \Uc \bi \Fc$  are isomorphic as Hopf algebras.
 \end{lemma}
 \begin{proof}
 Using the fact that the algebra structure of $U\bi F$ is the crossed product algebra,
  we first show that $X\rt 1$, $Y\rt 1$ and $1\rt \delta_k$ satisfy the identities (\ref{pres}).
 Indeed,
 \begin{align*}
& [Y\rt 1, X\rt 1]= YX\rt 1-XY\rt 1=X\rt 1, \\
 &\\
 &[X\rt 1, 1\rt \delta_k]= (X\rt1)(1\rt \delta_k)-(1\rt \delta_k)(X\rt 1)=\\
 &X\rt \delta_k-X\rt \delta_k-1\rt \delta_k\lt X=1\rt \delta_{k+1}, \\
 &\\
 &[Y\rt 1, 1\rt \delta_k]=(Y\rt1)(1\rt \delta_k)-(1\rt \delta_k)(Y\rt 1)=\\
 &Y\rt \delta_k-Y\rt \delta_k-1\rt \delta_k\lt Y=1\rt k\delta_k, \\
 &\\
 &[1\rt \delta_j,1\rt \delta_k]=1\rt \delta_j\delta_k-1\rt\delta_k\delta_j=0.
 \end{align*}
 Now we use the coalgebra structure of $ \Uc \bi \Fc$ which is of a
  cocrossed product coalgebra to check
   the identities (\ref{copr1})--(\ref{coun}):
 \begin{align*}
 &\Delta(X\rt 1)=X\ps{1}\rt X\ps{2}\ns{-1}\ot X\ps{2}\ns{0}\rt 1=\\
 &X\rt 1\ot 1\rt 1+1\rt X\ns{-1}\rt X\ns{0}\rt 1=\\
 &X\rt 1\ot 1\rt 1+1\rt X\ot 1\rt 1+1\rt \delta_i\ot Y\rt 1, \\
 &\\
  &\Delta(Y\rt 1)=Y\ps{1}\rt Y\ps{2}\ns{-1}\ot Y\ps{2}\ns{0}\rt 1=\\
   &Y\rt 1\ot 1\rt 1+1\rt Y\ns{-1}\rt Y\ns{0}\rt 1=\\
   &Y\rt1\ot 1\rt 1+1\rt Y\ot 1\rt 1, \\
   &\\
   &\Delta(1\rt \delta_1)=1\rt\delta_1\ps{1}\ot 1\rt \delta_1\ps{2}=1\rt \delta_1\ot 1\rt 1+1\rt 1\ot 1\rt \delta_1.
 \end{align*}
 Finally, using that the  counit of $\Uc \bi \Fc$ is
 $$\epsilon(h\rt f)=\epsilon(h)\epsilon(f) \, ,
 $$
 one easily verifies that
 \begin{equation*}
 \epsilon(X\rt1)=\epsilon(Y\rt 1)=\epsilon(1\rt \delta_k)=0.
 \end{equation*}
 \end{proof}

The maps  $\Psi$ and $\Psi^{-1}$  defined in \S 2.2, although
not shown to be cyclic (comp. Proposition \ref{map}) in the case at hand ,
are easily seen to preserve the cosimplicial
 structure. As a matter of fact the latter property holds
for any Hopf algebra whose underlying coalgebra is
 a cocrossed product. Thus,
   at least at the cosimplicial level, the diagonal subcomplex is still quasi-isomorphic
  to the Hochschild complex associated to the
bicosimplicial complex whose component in bidegree $(p,q)$ is
$\Fc^{\ot p}\ot   \Uc^{\ot q}$. More precisely, we are
 interested  in the cohomology of  the total complex of the following  bicomplex:
 \begin{center}
$\begin{xy}
\xymatrix{ \vdots&  \vdots & \vdots & \\
\Uc^{\ot 2}\ar[r]^\hb\ar[u]^\dvb& \Fc\ot   \Uc^{\ot 2}\ar[r]^\hb\ar[u]^\dvb&
 \Fc^{\ot 2}\ot   \Uc^{\ot 2}\ar[r]^\hb\ar[u]^\dvb &\dots  \\
 \Uc \ar[r]^\hb\ar[u]^\dvb& \Fc\ot   \Uc \ar[r]^\hb\ar[u]^\dvb&
\Fc^{\ot 2}\ot   \Uc \ar[r]^\hb\ar[u]^\dvb &\dots \\
 \bbc \ar[r]^\hb\ar[u]^\dvb& \Fc \ar[r]^\hb\ar[u]^\dvb& \Fc^{\ot 2} \ar[r]^\hb\ar[u]^\dvb &\dots   }
\end{xy}$
\end{center}

In order to evaluate it, we define a filtration on the total complex as follows:
  \begin{equation*}
  F_i= \bigoplus_{p\ge i}\bigoplus_{q\ge 0} \Fc^{\ot p}\ot   \Uc^{\ot q} \, .
  \end{equation*}

 \begin{lemma}
The $E_1$ term of the above spectral sequence is as follows:
\begin{eqnarray*}
E_1^{p,0}&=&\Fc^{\ot p} \, , \qquad
E_1^{p,1}\cong (\Fc^{\ot p}\ot \Cb X)\ot ( \Fc^{\ot p}\ot \Cb Y) \, ,\\
E_1^{p,2}&\cong& \Fc^{\ot p}\ot \Cb X\wedge Y \, , \qquad
E_1^{p,q} = 0 \, , \quad  \fl \,  q\ge 3 \, .
   \end{eqnarray*}
  \end {lemma}

  \begin{proof}
  The term $E_1$ is the cohomology of the columns with respect to the boundary $\vb$\,
   defined in (\ref{b}). In other words, $E_1^{p,q}=H^q(\Uc, \Fc^{\ot p})$ and since
 $\Uc$ coacts trivially on $\Fc^{\ot p}$, one gets
 $$E_1^{p,q} \, = \, H^q(\Uc,\bbc)\ot \Fc^{\ot p} \, .
 $$
It is well-known that (cf. e.g. \cite{ka, cm2}) the
 antisymmetrisation map induces an isomorphism
 between coalgebra Hochschild cohomology of the universal enveloping algebra  of a Lie algebra
 $\Fg$ and its exterior algebra $\wedge^\ast \Fg$.
 Accordingly, $E_1$ is given by: 
 
\begin{center}
$\begin{xy}
\xymatrix{ \vdots&\vdots& \vdots\\ 0&0&0\\
\wedge ^2\Fg \ar[u]& \Fc\ot \wedge ^2\Fg\ar[u]&
 \Fc^{\ot 2}\ot\wedge ^2\Fg\ar[u] &\dots  \\
\Fg \ar[u]^0& \Fc\ot \Fg \ar[u]^0&
\Fc^{\ot 2}\ot \Fg \ar[u]^0 &\dots \\ 
  \bbc \ar[u]^0 & \Fc \ar[u]^0 & \Fc^{\ot 2} \ar[u]^0 &\dots   }
\end{xy}$
\end{center}
 \end{proof}

\medskip

Since all maps are weight preserving our spectral sequence is graded,
\begin{equation*}
E_r=\bigoplus_{k\ge 0} E_r[k], \quad E_r[k]=\{\td f\ot   \td u\in
E_r \mid \quad \nm{\td f} + \nm{\td u}=k\}.
\end{equation*}

 Taking into consideration
 that we are working in the normalized complex, one sees that $E_1[1]$ 
 reduces to
 
\begin{center}
$\begin{xy}
\xymatrix{&&&\vdots&\vdots&\vdots & \\
&&&0&0&0& \dots\\
&&& \bbc X\wedge Y\ar[u]& 0\ar[u] &0\ar[u] & \dots  \\
E_1[1] &&& \bbc X\ar[u]^0&\bbc\delta_1\ot \bbc Y \ar[u]& 0\ar[u]&\dots \\
&&& 0 \ar[u]& \bbc \delta_1\ar[u]^0&0\ar[u]&\dots  }
\end{xy}$
\end{center}

and it then easily follows that $E_2[1]$ looks as follows :
\begin{center}
$\begin{xy}
\xymatrix{&&&\vdots&\vdots&\vdots & \\
&&&0&0&0& \dots\\
&&& \bbc X \wedge Y \ar[r]& 0\ar[r] &0\ar[r] & \dots  \\
E_2[1] &&& \bbc X\ar[r]^{\hb}\ar[r]&\bbc \delta_1\ot \bbc Y \ar[r]& 0\ar[r]&\dots \\
&&& 0 \ar[r]& \bbc \delta_1\ar[r]&0\ar[r]&\dots  }
\end{xy}$
\end{center}

Furthermore, since $\, \hb(X)= \delta_1\ot Y$, the next stage 
simplifies it further to only  two nontrivial classes:
\begin{center}
$\begin{xy}
\xymatrix{&&&\vdots&\vdots&\vdots &\vdots & \\
&&&0\ar[rrd]&0\ar[rrd]&0& 0& \dots\\
&&& \bbc X\wedge Y \ar[rrd]& 0\ar[rrd] &0&0 &\dots  \\
E_3[1] &&& 0 \ar[rrd]&0\ar[rrd]& 0&0 &\dots \\
&&& 0 & \bbc \delta_1&0& 0&\dots  }
\end{xy}$
\end{center}
\bigskip

\begin{proposition}\label{Hoch}
The weight $1$ component of the Hochschild cohomology of $H^\ast(\Hc_1)$
is generated by the classes
 $[\delta_1]$ and  $[X\ot Y-Y\ot X-Y\delta_1\ot Y]$.
\end{proposition}

\begin{proof}
From the above picture, one sees that in order to prove the statement
it suffices to convert  the nontrivial weight $1$
classes $[\delta_1]\in E_2^{1,0}$ and
 $X \wedge Y \in E_2^{0,2}$ into
  Hochschild  classes of $\Hc_1$.  This is a two-step process.
 First one needs to pass from the total complex to the
 diagonal of the bicosimplicial complex $\Fc^* \ot \Uc^*$,  via the
 Alexander-Whitney  map \eqref{AW}. The
 second move is the
 return to the Hochschild complex of
$\Hc_1$, via the map $\Psi^{-1}$ defined in  Proposition \ref{map}.

The result of performing the first move is the following:
 \begin{align*}
  &AW(X \ot Y)=\hd_2\hd_1(X \ot Y) = \hd_2(X\ps{-1}Y\ps{-1}\ot   X\ns{0}\ot Y\ns{0})=\\
  &\hd_2(1\ot   X\ot Y)+\hd_2(\delta_1\ot   Y \ot Y) =\\
 & 1\ot X\ns{-1}Y\ns{-1}\ot   X\ns{0}\ot Y\ns{0}+ \delta_1\ot Y\ns{-1}Y\ns{-1}\ot Y\ns{0}\ot Y\ns{0}=\\
 &1\ot \delta_1\ot Y\ot Y+ 1\ot 1\ot X\ot Y +\delta_1\ot 1\ot Y\ot Y ,
  \end{align*}
   and similarly
  $$AW(Y\ot X)=\delta_1\ot 1\ot Y\ot Y +1\ot \delta_1\ot  Y\ot Y+ 1\ot 1\ot  Y\ot X.
  $$
 This implies that
 $$AW(X\wedge Y)= 1\ot 1\ot X\ot Y-1\ot 1\ot Y\ot X.
 $$
  On the other hand
  $$\Psi^{-1}(1\ot 1\ot X\ot Y)=X\rt 1\ot Y\rt 1-Y\rt \delta_1
\ot Y\rt 1-Y\rt 1\ot Y\rt\delta_1, $$
    $$\Psi^{-1}(1\ot 1\ot Y\ot X)=Y\rt 1\ot X\rt 1-Y\rt 1
\ot Y\rt \delta_1 \, ,$$
 and after the usual identification we obtain
$$\Psi ^{-1}(1\ot1\ot X\wedge Y)= X\ot Y- Y\ot X -Y\delta_1\ot Y.
 $$
 In the case of $[\delta_1]\in E_2^{1,0}$
 one has $AW(\delta_1)=-\delta_1\ot 1$,
 and hence $\Psi^{-1}(\delta_1\ot 1)=1\rt\delta_1$,
 identified to $\delta_1$.
\end{proof}
\medskip

We are now in a position to give an intrinsic proof of one of the key results
in \cite{cm2}, without
appealing to Gelfand-Fuks cohomology.
\medskip

\begin{theorem}\label{HP(H_1)}
The periodic cyclic cohomology of the Hopf algebra $\Hc_1$
with respect to the MPI $(\delta,1)$ is the following:
$$HP_{(\delta,1)}^{\odd}(\Hc_1)= \Cb\cdot  [\delta_1] \, , \qquad
HP^{\ev}_{(\delta,1)}(\Hc_1)= \Cb\cdot  [X\ot Y-Y\ot X-\delta_1Y\ot Y]\, .$$
\end{theorem}

\begin{proof}
From Corollary \ref{percartan} we know that $HP^\ast_{(\delta,\sigma)}(\Hc_1)$ 
can be computed solely from the weight $1$
mixed complex $({\Hc_1}_\nr [1], b,B)$.  By the above proposition
  the Hochschild cohomology of the latter complex is generated by
  the Hochschild classes of $\delta_1$
 and $X\ot Y-Y\ot X-Y\delta_1\ot Y$. As
 both of them are cyclic, the Connes boundary
$B$ is zero on the Hochschild cohomology.

 To obtain the stated form of the result it remains to note that the classes
 $[X\ot Y-Y\ot X-\delta_1Y\ot Y ]$ and $ [X\ot Y- Y\ot X-Y\delta_1\ot Y]$
coincide, since 
$$
Y\delta_1\ot Y - \delta_1Y\ot Y \, = \, \d_1 \ot Y \, = \, b(X)
\, = \, (b + B)(X) \, .
$$
\end{proof}
\medskip

A similar result holds for the quotient Hopf algebra
$\Hc_{1 s} $ = $\Hc_1/\Sc$, where $\Sc$ is the ideal generated by the `Schwarzian'
$\d_2'=\delta_2-\frac{1}{2}\d_1^2$. The Hopf algebra $\Hc_{1 s} $ plays the
crucial role in the extension of the Rankin-Cohen brackets to modular
Hecke algebras (cf.~\cite{cm5'}). 
\smallskip

Evidently, as an algebra,
 $\Hc_{1 s} $ is generated by (the classes \textit{mod} $\Sc$ of)
 $X$, $Y$ and $Z$=$ \d_1$,
 subject to the relations
\begin{equation*}
[Y,X]=X, \quad [Y,Z]=Z, \quad [X,Z]=\frac{1}{2}Z^2; 
\end{equation*}
its coalgebra structure is defined by 
\begin{align*}
&\D(X)=1\ot X+X\ot 1 +Z\ot Y,\\
&\D(Y)=Y\ot 1+ 1\ot Y,\\
&\D(Z)=Z\ot 1+1\ot Z, \\
&\epsilon(X)=\epsilon(Y)=\epsilon(Z)=0 \, ,
\end{align*}
and the antipode is defined by 
\begin{align*}
S(X)= -X+ZY, \quad S(Y)=-Y, \quad S(Z)=-Z.  
\end{align*}
\medskip

\begin{theorem}\label{HP(H_s)}
The periodic cyclic cohomology of the Hopf algebra $\Hc_{1 s}$
with respect to the MPI $(\delta,1)$ is the following:
$$HP_{(\delta,1)}^{\odd}(\Hc_{1 s})= \Cb \cdot [Z] \, , \qquad
HP^{\ev}_{(\delta,1)}(\Hc_{1 s})= \Cb \cdot  [X\ot Y-Y\ot X-ZY\ot Y]\, .$$
\end{theorem}
\smallskip

\begin{proof}  We first reconstruct  $\Hc_{1 s} $ as a bicrossed product Hopf algebra,
emulating Lemma \ref{bicp}.  Indeed, the algebra 
  $\Uc=\Uc({\Fg}_{-})$ acts from the right on the polynomial algebra $\bbc[Z]$  via
\begin{equation*}
Z\lt X=-\frac{1}{2}Z^2, \quad Z\lt Y=-Z \, ,
\end{equation*}
and turns it into a right $\Uc$-module algebra. 
On the other hand, $\bbc[Z]$, viewed as the free Hopf algebra generated by the primitive
element $Z$,
 coacts on $\Uc$ via $\rho: \Uc\ra \bbc[Z]\ot \Uc$ defined by
\begin{equation*}
\rho(X)=1\ot X+Z\ot Y, \quad \rho(Y)=1\ot Y \, ,
\end{equation*}
Endowed with this coaction $\Uc$ is a 
$\bbc[Z]$-comodule coalgebra. Moreover, with respect to
the above action and coaction,
$\Uc$ and $\bbc[Z]$ form 
a matched pair of Hopf algebras (cf.~\cite[page 322]{maj}). 
One can then form  the bicrossed product Hopf algebra $\Uc\bi \bbc[Z]$, which is clearly
isomorphic to   $\Hc_{1 s}$, via the isomorphism
identifying $X\rt 1$, $Y\rt 1$, and $1\rt Z$ with $X$, $Y$, and $Z$ respectively.
 \smallskip
 
To obtain the periodic cohomology $HP_{(\delta,1)}^{\ast}(\Hc_{1 s})$, it suffices to
repeat verbatim the proof of Theorem \ref{HP(H_1)}, after
replacing $\Fc$ by $\bbc[Z]$.
\end{proof}

\subsection{Hopf cyclic cohomology of $\Hc_1^\dag$} \label{S4.1}

In this  section we compute the periodic cyclic cohomology of the Hopf algebra
$\Hc_1^\dag$ defined in $\S$\ref{S1.2}, with coefficients in
$\Kc= \bbc[\sigma,\sigma^{-1}] $ viewed as an SAYD module.
\medskip

In order to make the method of $\S$\ref{S2} applicable, we shall represent
$\Hc_1^\dag$ as a cocrossed product coalgebra. To this end, we define a
 coaction of $\Kc$ on $\Hc_1$,
\begin{equation}\label{comodule}
\rho: \Hc\ra \Kc\ot \Hc, \quad \rho(h)=\sigma^{\nm{h}}\ot h.
\end{equation}

\begin{lemma}\label{K-comodule}
The map $\rho: \Hc\ra \Kc\ot \Hc$
defines a coaction, which makes $\Hc_1$ a $\Kc$-comodule Hopf algebra.
\end{lemma}
\begin{proof}
The first claim is obvious. Next,
$\Hc_1$ is a $\Kc$-comodule algebra  because
$\rho(1)=1\ot 1$ and
$$\rho(hg)=\sigma^{\nm{hg}\ot hg}=\sigma^{\nm{h}+\nm{g}}\ot
hg=(\sigma^{\nm{h}}\ot h)(\sigma^{\nm{g}}\ot g)=\rho(h)\rho(g).$$

To prove that it is also a comodule coalgebra, we note that
$$h\ns{-1}\ot\epsilon(h\ns{0})=\sigma^{\nm{h}}\ot \epsilon(h)=1\ot
\epsilon(h),
$$
where in the latter equality we used the fact that
$\epsilon(h)=0$ if $\nm{h}>0$.

Finally, the coaction and the coproduct do commute:
 \begin{align*}
&  h\ps{1}\ns{-1} h\ps{2}\ns{-1}\ot h\ps{1}\ns{0}\ot h\ps{2}\ns{0}
=\sigma^{\nm{h\ps{1}}}\sigma^{\nm{h\ps{2}}}\ot h\ps{1}\ot h\ps{2}\\
&=\sigma^{\nm{h\ps{1}}+\nm{h\ps{2}}}\ot  h\ps{1}\ot
h\ps{2}=\sigma^{\nm{h}}\ot h\ps{1}\ot h\ps{2}=h\ns{-1}\ot
h\ns{0}\ps{1}\ot h\ns{0}\ps{2}.
\end{align*}
\end{proof}

Since $\Hc_1$, $\Kc$  and $\rho$ satisfy the condition of Lemma
\ref{crossed}, one can form the Hopf algebra $\Hc_1\rt \Kc$.

\begin{lemma}
The Hopf algebras $\Hc_1^\dag$ and $\Hc_1\rt \Kc$ are isomorphic.
\end{lemma}
\begin{proof}
The map $\Phi: \Hc_1^\dag\ra \Hc_1\rt\Kc$, defined on
monomials  by $\Phi(\sigma^mh)=h\rt \sigma^m$, extends by
multiplicativity to an algebra map. We check that
it commutes with coproducts
and counits. Indeed,
\begin{align*}
&\Delta_{\Hc_1\rt \Kc}(\Psi(\sigma^mh))=
\Delta_{\Hc_1\rt\Kc}(h\rt\sigma^m)= && &&\\
&=h\ps{1}\rt h\ps{2}\ns{-1}\sigma^m\ps{1}\ot h\ps{2}\ns{0}\rt\sigma^m\ps{2}\\
&= h\ps{1}\rt \sigma^{\nm{h\ps{2}}+m}\ot h\ps{2}\rt \sigma^m=\\
&=\Psi\ot \Psi(\sigma^{m+\nm{h\ps{2}}}h\ps{1}\ot
\sigma^mh\ps{2})=\Psi\ot\Psi(\Delta_{\Hc_1^\dag}(\sigma^mh))
\end{align*}
and
\begin{align*}
 &\epsilon_{\Hc_1\rt\Kc}(\Psi(\sigma^mh))
=\epsilon_{\Hc_1\rt\Kc}(h\rt\sigma^m)=
\epsilon_{\Hc_1}(h)\epsilon_\Kc(\sigma^m)=\epsilon_{\Hc_1^\dag}(h).
\end{align*}
\end{proof}
\medskip

We recall that $(\delta,1)$ is a MPI for $\Hc_1$, with
 $\delta(X)=\delta(\delta_k)=0$ and $\delta(Y)=1$. On the
other hand, the MPIs for $\Kc$ fall into two classes: those of the form
$(\beta, 1)$,
where $\beta$ is any character of $\Kc$, and those of
the form $(\beta_{2k\pi
 i/n},\sigma^{n})$,   $k \in \mathbb{Z}$, $0\neq n\in \mathbb{Z}$, where
  $\beta_z(\sigma)=e^z$,  $z\in \bbc$.

 Among the MPIs listed above
 only those of the form $(\epsilon,\sigma^k)$, $k\in \mathbb{Z}$, for $\Kc$,
 together with $(\delta,1)$, for $\Hc_1$, satisfy the conditions of Proposition
 \ref{MPI}, and thus combine to give
 $(\delta\ot\epsilon, 1\ot \sigma^k)$,  $k\in \mathbb{Z}$, MPIs for
 $\Hc_1^\dag$.
\bigskip

 Let us compute the Hopf cyclic cohomology of $\Hc_1^\dag$
with coefficients in the MPI $(\delta,\sigma^k)$, where we identify
 $\delta \equiv \delta\ot\epsilon$  and $\sigma^k \equiv 1\rt \sigma^k $.
\smallskip

Since by Theorem \ref{mix} the mixed complex $(\Hc_1^\dag, b,B)$
is quasi-isomorphic to $(\Tot(\FC), b_T, B_T)$, we start by 
 computing the Hochschild
 cohomology of $\Tot(\FC)$.
 To this end,  we introduce the following filtration on  $\Tot(\FC)$:
\begin{equation*}
F^n=\bigoplus_{q\ge n}\bigoplus_{p\ge 0}\Kc^{\ot p}\ot \Hc_1^{\ot q}.
\end{equation*}
Thus, $E^{p,q}_1$ is the $p$th cohomology of the $q$th row in the diagram
below:
 \begin{center}
$\begin{xy}
\xymatrix{ &&&\vdots&  \vdots & \vdots & \\
&&&\Hc_1^{\ot 2}\ar[r]^\hb& \Kc\ot   \Hc_1^{\ot 2}\ar[r]^\hb&
 \Kc^{\ot 2}\ot   \Hc_1^{\ot 2}\ar[r]^\hb &\dots  \\
E_0&&& \Hc_1\ar[r]^\hb& \Kc\ot   \Hc_1\ar[r]^\hb&
\Kc^{\ot 2}\ot   \Hc_1\ar[r]^\hb &\dots \\
&&& \bbc \ar[r]^\hb& \Kc \ar[r]^\hb& \Kc^{\ot 2} \ar[r]^\hb&\dots   }
\end{xy}$
\end{center}

This cohomology  is precisely
the Hochschild cohomology of $\Kc$ with coefficients in $\Hc_1^{\ot q}$,
viewed as a left $\Kc$-comodule via the coaction
\begin{equation*}
\td\rho: \Hc_1^{q}\ra \Kc\ot \Hc_1^{\ot q}, \quad \td\rho(\td h)=\sigma^{k+\nm{\td h}}\ot \td{h}.
\end{equation*}

 To compute it,
 we use the
 free $\Kc$-comodule resolution for the trivial  $\Kc$ comodule $\bbc$ given in
 \cite{cr}, namely:
\begin{equation} \nonumber
0\rightarrow \bbc \overset{\eta}{\rightarrow} \bbc[\sigma,\sigma^{-1}]
\overset{\theta}{\rightarrow} \bbc[\sigma,\sigma^{-1}]
 \overset{\gamma}{\rightarrow} \bbc[\sigma,\sigma^{-1}]\overset{\theta}
{\rightarrow} \bbc[\sigma,\sigma^{-1}] \overset{\gamma}{\rightarrow} \dots
\end{equation}
where $\eta=1_\Kc$, while $\theta$ and $\gamma$ are defined on generators by
\begin{equation*}
\theta(1)=0,\quad  \theta(\sigma^i)=\sigma^i,\;\; i\neq 0,
\end{equation*}
\begin{equation*}
\gamma(1)=1,\quad \gamma(\sigma^j)=0,\,\,  j\neq 0.
\end{equation*}

One then finds that  $\cotor^\ast _\Kc(\bbc, \Hc_1^{\ot q})= \Hc_1^{\ot q}
 $ with the following coboundary:
$$\Hc_1^{\ot q} \overset{\td\theta}{\rightarrow} \Hc_1^{\ot q}
\overset{\td\gamma}{\rightarrow} \Hc_1^{\ot q}
 \overset{\td \theta }{\rightarrow} \Hc_1^{\ot q}\overset{\td\gamma}{\rightarrow}  \dots, $$
where  $\td \theta$ and $\td \gamma$ are defined by
\begin{equation*}
\td\theta(\td h)=0,  \, \, \text{if}\,\, \,\,\nm{\td h}=-k,
  \qquad  \td\theta(\td h)=\td h,\,\,\text{if }\,\, \nm{\td h}\neq-k,
\end{equation*}
\begin{equation*}
\td\gamma(\td h)=\td h, \,\,\text{if }\,\, \nm{\td h}=-k\quad
\td\gamma(\td h)=0,\,\,\,\text{if}\,\,\, \nm{\td h}\neq -k.
\end{equation*}
So the resulting cohomology is as follows:
\begin{equation*}
E^{p,q}_1=0, \, p>0, \quad E_1^{0,q}={\Hc^q_1}_\nr[-k]:=\{\td h\in\Hc_1^{\ot q}\mid \nm{\td h}=-k \} \, ,
\end{equation*}
or schematically
\begin{center}
$\begin{xy}
\xymatrix{ &&&\vdots &\vdots  \\
&&&\Hc_1^{\ot 2}[-k]&0&\dots \\
E_1&&& \Hc_1[-k]&0&\dots\\
&&& 0&0&\dots}
\end{xy}$
\end{center}

To compute the $E_2$ term of this spectral sequence one notes that the vertical
 operators are exactly those of the
 cocyclic submodule ${\Hc_1}_\nr[-k]$,

 \begin{center}
$\begin{xy}
\xymatrix{ &&&\vdots  &\vdots &\\
&&&\Hc_1^{\ot 2}[-k] \ar[u]^{b_{\Hc_1}} &0&\dots\\
E_2&&& \Hc_1[-k]\ar[u]^{b_{\Hc_1}}&0&\dots\\
&&& \bbc[-k]\ar[u]^{b_{\Hc_1}} &0&\dots}
\end{xy}$
\end{center}
Therefore
\begin{equation*}
E_2^{p,q}=0, \,\, p\ge 0,\qquad E_2^{0,q} =H^q_{(\delta,1)}({\Hc_1}_\nr[-k])\, ,
\end{equation*}
where the lower subscript for the cohomology groups stands for the coefficients.
\smallskip

Thus, the Hochschild cohomology of $\Hc_1^\dag$
 with coefficients in  $(\delta,\sigma^k)$ is  
\begin{equation*}
H^n_{(\delta,\sigma^k)}(\Hc^\dag_1)= H^n_{(\delta,1)} ({\Hc_1}_\nr[-k]).
\end{equation*}
Since
the vertical operators coincide with those belonging to
 the cocyclic structure for $(\Hc_1; \delta, 1)$ we conclude that
\begin{equation}\label{HP}
HP^\ast_{(\delta,\sigma^k)}(\Hc^\dag_1)=HP^\ast_{(\delta,1)}({\Hc_1}_\nr[-k])
\end{equation}
\medskip

Remark that all the preceding arguments work as well when
the $\Hc_1^\dag$ is replaced by any of the  `finite' cyclic covers
$\Hc_1^{\dag|N}$, with $N > 1$.
\medskip
 
\begin{theorem}\label{final1}
$1^0$. Let $N>1$. Of all the MPIs $(\delta,\sigma^k)$, $k\in \mathbb{Z}$,  only
$(\delta,\sigma^{-1})$ yields nontrivial periodic cyclic
cohomology:
$$HP^*_{(\delta,\sigma^k)}(\Hc_1^{\dag|N})=0 \quad {\rm if} \quad
k\neq -1, \quad {\rm and} \quad
HP^*_{(\delta,\sigma^{-1})}(\Hc_1^{\dag|N})\cong
HP^*_{(\delta,1)}(\Hc_1) .
$$
$2^0$. The periodic cyclic cohomology of $\Hc_1^{\dag|N}$,  $N > 1$,
 with coefficients in
$(\delta,\sigma^{-1})$ is generated by  $[TF^\dag]$  in the even degree  and
by $[\delta_1^\dag]$ in the  odd degree, where
\begin{equation*}
TF^\dag=\sigma^{-1} X\ot \sigma^{-1} Y-Y\ot \sigma^{-1} X-\sigma^{-1}
\delta_1Y\ot  \s^{-1}Y, \quad \delta_1^\dag=  -\sigma^{-1}\delta_1.
\end{equation*}
\end{theorem}

\begin{proof}
For the first claim, it  suffices to apply Corollary \ref{percartan} to
the right  hand side of (\ref{HP}).

To prove the second, we start from the known generators
  $[\delta_1]$ and $[TF]= [X\ot Y-Y\ot X-\delta_1Y\ot Y]$
for the periodic cyclic cohomology of $\Hc_1$ and go through
the two moves required by the above routine to convert them into
periodic cyclic cohomology  classes for $\Hc^\dag_1$. The first
is to transfer back these classes to the total complex
via Alexander-Whitney  map \eqref{AW}, and
the second amounts to applying $\Psi^{-1}$.

Concerning the first move, we note that
for $\td h\in E^{0,q}$ one has
\begin{equation*}\hd_2\hd_1(\td h)=  \sigma^{(\nm{\td h}-1)}\ot \sigma^{(\nm{\td h}-1)}\ot   \td h \, ;
\end{equation*}
indeed,
\begin{align*}
&\hd_2\hd_1(h^1\ot h^2)=\hd_2(\sigma^{-1}h^1\ns{-1}h^2\ns{-1}\ot   h^1\ns{0}\ot h^2\ns{0})=\\
&\hd_2(\sigma^{\nm{h^1}+\nm{h^2}-1}\ot   h^1\ot h^2)= \sigma^{\nm{h^1}+\nm{h^2}-1}
\ot \sigma^{-1}h^1\ns{-1}h^2\ns{-1}\ot   h^1\ns{0}\ot h^2\ns{0}\\
&=\sigma^{\nm{h^1}+\nm{h^2}-1}\ot\sigma^{\nm{h^1}+\nm{h^2}-1}\ot   h^1\ot h^2.
\end{align*}

Applying this to $TF$ one obtains
\begin{equation*}
AW(TF)=1\ot 1 \ot   X\ot Y-1\ot 1 \ot   Y\ot X-1\ot 1 \ot   \delta_1Y\ot Y
\end{equation*}
In the second move, we get
\begin{align*}
&\Psi^{-1}(1\ot 1\ot   X\ot Y)= X\rt\sigma^{-1}\ot Y\rt \sigma^{-1}=\sigma^{-1} X\ot \sigma^{-1} Y,\\
&\Psi^{-1}(1\ot 1\ot   Y\ot X)= Y\rt 1\ot X\rt \s^{-1}=Y\ot \sigma^{-1} X,\\
&\Psi^{-1}(1\ot 1\ot   \delta_1Y\ot Y)= \delta_1Y\rt \sigma^{-1}\ot Y\rt \s^{-1}=
\sigma^{-1}\delta_1Y\ot  \s^{-1}Y,\\
\end{align*}
which implies that
 $$\Psi^{-1} \circ AW(TF)=\sigma^{-1} X\ot \sigma^{-1} Y-Y\ot \sigma^{-1} X-\
 \sigma^{-1}\delta_1Y\ot \sigma^{-1}Y.
 $$

Repeating the procedure for the class $\delta_1$,one obtains
$$AW(\delta_1)=-1\ot \delta_1 \quad \text{and}  \quad\Psi^{-1}(1\ot \delta_1)=
\delta_1\rt\sigma^{-1}=\sigma^{-1}\delta_1 \, ,
$$
and therefore
 $$ \Psi^{-1}\circ AW(\delta_1)=-\sigma^{-1}\delta_1\, .
 $$
\end{proof}


\subsection{Hopf cyclic cohomology of $\Hc_{CK}$ and $\Hc_{CK}^\dag$} \label{S4.3}

In this subsection we apply the same method as in the previous two
in order
to calculate the Hopf cyclic cohomology of the extended
Connes-Kreimer Hopf algebra $\Hc_{CK}$.
\medskip

We begin by recalling from \cite{ck} the basic ingredients from which
the Hopf algebra $\Hc_{CK}$ is constructed.

A \textit{rooted tree} $T$ is a finite, connected and simply connected
 $1$-dimensional simpicial complex with a base point (root)
$\ast\in \Delta^0=\{\text{set of vertices of}\,\, T\}$.
The set of rooted trees, up to isomorphism, is denoted by $\Sigma$.
It is graded and by weight, which is defined as the integer
\begin{equation*}
\nm{T}= Card \0D(T)=\# \,\text{of vertices of }\, T \, .
\end{equation*}
A \textit{cut} of a rooted tree $T$  is a subset $c\subset \1D(T)$ of the set
of edges of $T$.  A cut is called {\it simple} if for any $x\in
\0D(T)$ the  path $(\ast,x)$ only contains one element of $c$. To
any simple cut of a rooted tree $T$ one assigns  a trunk denoted by
$R_c(T)$ and a collection of cut branches denoted by $P_c(T)$.

Let $\Hc_{rt}$ be the free commutative algebra generated by the
symbols
\begin{equation} \nonumber
\{ \delta_T\, ; \quad T\in \Sigma \, \} \, .
\end{equation}
One defines a coproduct on $\Hc_{rt}$ by extending to an algebra homomorphism
the map defined on the generators of $\Hc_{rt}$  by
\begin{equation*}
\Delta(\delta_T)=\delta_T\ot 1+ \sum_{\text{simple
cuts}}\left(\prod_{P_c(T)}\delta_{T_i}\right)\ot \delta_{R_c(T)}.
\end{equation*}
\medskip

We let $\Fg$,  $X$, and $Y$ be the same as in (\ref{XY}), and let
$\Uc:=\Uc(\Fg)$ denote the universal enveloping algebra of the Lie algebra
$\Fg$. It acts from right on $\Hc_{rt}$ as follows:
\begin{align*}
\d_T\lt Y= -\nm{T}\d_T, \\
\d_T\lt X=- N\d_T:=-\sum \d_{T'} \, ,
\end{align*}
where the trees $T'$ are obtained by adding one vertex and one edge to $T$ in all possible ways without changing the base point. One then extends the action of $\Uc$
on $\Hc_{rt}$, so that $\Hc_{rt}$ becomes a $\Uc$-module algebra.

On the other hand, as in the case of $\Hc_1$,
one defines a left coaction of $\Hc_{rt}$ on $\Uc$,
$\rho: \Uc\ra \Hc_{rt}\ot \Uc$, by setting
\begin{align*}
&\rho(X)= 1\ot X + \delta_\ast\ot Y \, , \\
&\rho(Y)=1\ot Y \, ,
\end{align*}
and then extending the map $\rho$ to $\Uc$ via the rule (\ref{bi}).

Equipped with this coaction $\Uc$ is a left $\Hc_{rt}$-comodule
algebra. Moreover the action and coaction thus defined
fulfill  the axioms of a matched pair~\cite[p. 232]{maj}.

   One then forms a bicrossed Hopf algebra $\Hc_{CK} \equiv \Uc\bi \Hc_{rt}$,
   whose algebra structure is the crossed product algebra $ \Uc\ltimes  \Hc_{rt}$,
   and coalgebra structure is $\Uc\rt \Hc_{rt}$, with antipode
\begin{equation*}
S(u\rt h)=(1\rt S_{\Hc_{rt}}(u\ns{-1}h))(S_{\Uc}(h\ns{0})\rt 1).
\end{equation*}

We identify $X\rt 1$, $Y\rt 1$ and $1\rt\d_T$ with $X$, $Y$, and $\d_T$ respectively.
One can easily see that, with the character $\d$ defined as before,
\begin{equation*}
\delta(X)=0, \quad \d(Y)=1, \quad \d(\d_T)=0 \, ,
\end{equation*}
 $(\delta,1)$ is a MPI for $\Hc_{CK}$.

 So the Hopf cyclic cohomology of $\Hc_{CK}$ with coefficients  in $(\d,1)$ is well defined,
 and  to compute it we
only need to repeat the method of the preceding subsections.
 \medskip

Let  $\Hc_{CK\nr} $ be the cocyclic module associated to $(\Hc_{rt} ; \d,1)$.
Since the operators defining the cyclic structure are weight preserving,
 the cocyclic module is
graded decomposes into eigenmodules
\begin{equation*}
\Hc_{CK\nr}= \bigoplus_{k\ge 0}\Hc_{CK\nr} [k] \, ,
\end{equation*}
where $\Hc_{CK\nr} [k]$ consists of all elements in $\Hc_{CK\nr}[k]$ of weight $k$.

Furthermore, as  $Y$ is a primitive  element of $\Hc_{CK}$ and
since  $\ad_Y$ implements the grading,
using the the Cartan homotopy formula  \ref{percartan},  one obtains that
  \begin{equation*}
 HP^\ast(\Hc_{CK\nr}[1])=HP^\ast_{(\d,1)}(\Hc_{CK}) \, .
  \end{equation*}

 We first  compute the Hochschild cohomology of left hand side.
  To this end we again use the bicosimplicial module
 $$\FC^{p,q}= \Hc_{rt}^{\ot p}\ot   \Uc^{\ot q}, $$
 with the cosimplicial structure defined by (\ref{hd0})--(\ref{vs}).

The following defines a  filtration on the total complex of the bicomplex $(\FC , \hb, \vb)$,
where the horizontal boundary  $\hb$ and and the vertical boundary
 $\vb$ are defined in (\ref{b}):
\begin{equation*}
 F_i= \bigoplus_{p\ge i}\bigoplus_{q\ge 0} \Hc_{rt}^{\ot p}\ot   \Uc^{\ot q} \,
\end{equation*}
and the associated spectral sequence $E_r$
converges to the Hochschild cohomology of the total complex
 $\Tot(\FC)$.

 The spectral sequence $E_r$ inherits the grading of $\Hc_{CK}$, since
all cosimplicial  maps,  as well as $\Psi $ and $\Psi^{-1}$ (defined in \S 2.2),
 are weight preserving.  As before, in view of the Cartan homotopy formula,
 it suffices to focus on its weight $1$ component $E_r[1]$.

    Using again the fact that the
     antisymmetrisation map induces an isomorphism
 between coalgebra Hochschild cohomology of the universal enveloping algebra  of a Lie algebra
 $\Fg$ and its exterior algebra $\wedge^\ast \Fg$, one obtains the term $E_1$ of the
 spectral sequence as follows.

\begin{proposition}
  $E_1^{p,0}=\Hc_{rt}^{\ot p}$, $E_1^{p,1}=
  (\Hc_{rt}^{\ot p}\ot X)\ot ( \Hc_{rt}^{\ot p}\ot Y)$,
   $ E_1^{p,q}=0$ for $q\ge 3$ and $E_1^{p,2}= \Hc_{rt}^{\ot p}\ot X\wedge Y$.
  \end {proposition}

Since up to a constant factor $\delta_\ast$ is the only element of weight $1$ 
in $\Hc_{rt}$,
 and taking into consideration that
we are always working in the normalized complex,
we obtain the following picture for $E_2[1]$:

  \begin{center}
$\begin{xy}
\xymatrix{ \bbc X\wedge Y \ar[r]& 0\ar[r] &0\ar[r] & \dots  \\
  \bbc[X]\ar[r]^{\hb}\ar[r]&\Cb \delta_\ast \ot \bbc Y \ar[r]& 0\ar[r]&\dots \\
 0 \ar[r]& \bbc \delta_\ast \ar[r]&0\ar[r]&\dots  }
\end{xy}$
\end{center}
Because $\hb(X)=\delta_\ast\ot Y$, one sees
 that only $[\delta_\ast]$, and $X\wedge Y$  survive
 in $E_2[1]$.
By the same arguments as in the
 proof of Proposition \ref{Hoch}, we obtain that the
   Hochschild cohomology in weight  $1$ is generated by the classes
 $$[\delta_\ast] \quad {\rm  and}  \quad [X\ot Y-Y\ot X-Y\delta_\ast\ot
Y] \, .
$$
Finally, one computes the periodic cyclic cohomology of $ \Hc_{CK}$
exactly as in the case of $\Hc_1$. In particular, the
proof of the following result is identical to that of
 Theorem \ref{HP(H_1)}.

\begin{theorem} \label{final2}
The periodic cyclic cohomology of the Hopf algebra $\Hc_{CK}$
with respect to the MPI $(\delta,1)$ is the following:
$$HP_{(\delta,1)}^1(\Hc_{CK})=\bbc \cdot [\delta_\ast]~~\text{and}~~
HP^0_{(\delta,1)}(\Hc_{CK})=\bbc \cdot [X\ot Y-Y\ot X-\delta_\ast Y\ot Y].$$
\end{theorem}
\medskip

We end with the observation that the Hopf algebra  $\Hc_{CK}$,
admits its own `cyclic coverings' too.
Indeed, as it was done for $\Hc_1$ in (\ref{comodule}),
one can
endow it with a similar coaction
$\rho: \Hc_{CK}\ra \Kc\ot \Hc_{CK}$ ,
\begin{equation*}
 \rho(h)= \sigma^{\nm{h}}\ot h \, .
\end{equation*}
One checks that this coaction turns  $\Hc_{CK}$ into a $\Kc$-comodule Hopf algebra,
as in the proof of Lemma \ref{K-comodule}.  We can thus form the cocrossed
product Hopf algebra
\begin{equation*}
\Hc_{CK}^\dag:=\Hc_{CK}\rt \Kc.
\end{equation*}
The periodic Hopf cyclic cohomology of $\Hc_{CK}^\dag$ as well as of the
finite cyclic covers  $\Hc_{CK}^{\dag|N}$, $N > 1$,
with coefficients
in the MPIs  $(\delta,\sigma^k)$, $k\in \mathbb{Z}$,
can be  computed in exactly the same way as for
 $\Hc_1^\dag$ (see Theorem \ref{final1}). We record below
 the analogous results.
 \medskip

\begin{theorem} \label{final3}
$1^0$. Let $N > 1$. Of all the MPIs $(\delta,\sigma^k)$, $k\in \mathbb{Z}$,  only
$(\delta,\sigma^{-1})$ yields nontrivial periodic cyclic
cohomology:
$$HP^*_{(\delta,\sigma^k)}(\Hc_{CK}^{\dag|N})=0 \quad {\rm if} \quad
k\neq -1, \quad {\rm and} \quad
HP^*_{(\delta,\sigma^{-1})}(\Hc_{CK}^{\dag|N})\cong
HP^*_{(\delta,1)}(\Hc_{CK}) .
$$
$2^0$. The periodic cyclic cohomology of $\Hc^{\dag|N}_{CK}$  with
coefficients in $(\delta,\sigma^{-1})$ is generated by  $[TF^\dag]$
in the even degree  and by $[\delta_\ast^\dag]$ in the  odd degree,
where
\begin{equation*}
TF^\dag=\sigma^{-1} X\ot \sigma^{-1} Y-Y\ot \sigma^{-1} X-\sigma^{-1}
\delta_\ast Y\ot \sigma^{-1} Y, \quad \delta_\ast^\dag=  -\sigma^{-1}\delta_\ast.
\end{equation*}
\end{theorem}
\bigskip

\end{document}